\UseRawInputEncoding

\documentclass{article}

\usepackage[numbers]{natbib}
\usepackage{graphicx}
\usepackage{notoccite} % To control the numbering of the references
\usepackage[ruled, vlined, lined, linesnumbered, commentsnumbered]{algorithm2e} % For the pseudo-code
\usepackage{multirow}
\newcommand{\rom}[1]{\lowercase\expandafter{\romannumeral #1\relax}}
\usepackage{amsmath,amssymb,amsfonts}
\usepackage{adjustbox}
\usepackage{mathtools}
\usepackage{caption}
\usepackage[figuresleft]{rotating}
\usepackage{academicons}
\usepackage{xcolor}

\DeclareCaptionLabelFormat{continued}{#1~#2 (Continued)}
\usepackage[margin=1.5in]{geometry}
\usepackage{authblk}
\usepackage{hyperref}

\newbox{\myorcidaffilbox}
\sbox{\myorcidaffilbox}{\large\includegraphics[height=1.7ex]{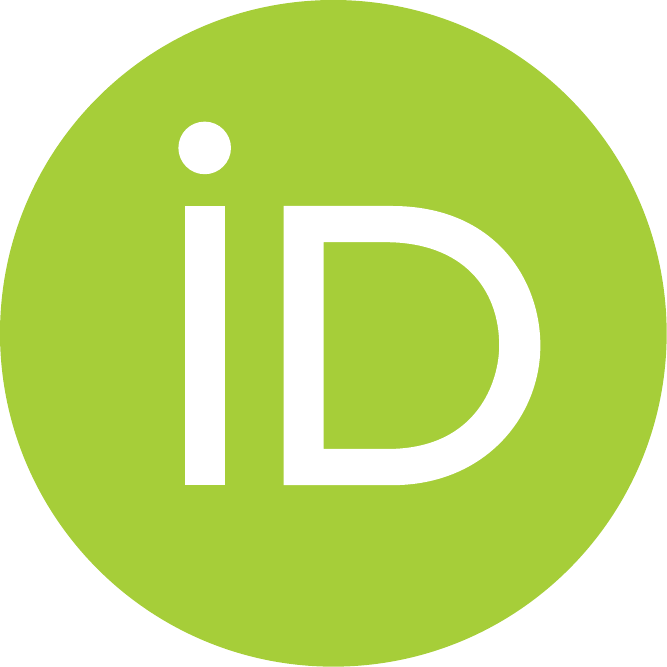}}
\newcommand{\orcidaffil}[1]{%
  \href{https://orcid.org/#1}{\usebox{\myorcidaffilbox}}}

\begin{document}

\title{%
Pareto-like Sequential Sampling Heuristic for Global Optimisation
\\~\\ \small [Preprint]}

\author[1]{Mahmoud Shaqfa \orcidaffil{0000-0002-0136-2391}}
\author[2]{Katrin Beyer \orcidaffil{0000-0002-6883-5157}}
\affil[1,2]{Earthquake Engineering and Structural Dynamics Laboratory (EESD), \'Ecole polytechnique f\'ed\'erale de Lausanne (EPFL), CH-1015 Lausanne, Switzerland.}

\date{}
\maketitle

% ********* Abstract and Keywords *********
\begin{abstract}
In this paper, we propose a simple global optimisation algorithm inspired by Pareto's principle. This algorithm samples most of its solutions within prominent search domains and is equipped with a self-adaptive mechanism to control the dynamic tightening of the prominent domains while the greediness of the algorithm increases over time (iterations). Unlike traditional metaheuristics, the proposed method has no direct mutation- or crossover-like operations. It depends solely on the sequential random sampling that can be used in diversification and intensification processes while keeping the information-flow between generations and the structural bias at a minimum. By using a simple topology, the algorithm avoids premature convergence by sampling new solutions every generation. A simple theoretical derivation revealed that the exploration of this approach is unbiased and the rate of the diversification is constant during the runtime. The trade-off balance between the diversification and the intensification is explained theoretically and experimentally. This proposed approach has been benchmarked against standard optimisation problems as well as a selected set of simple and complex engineering applications. We used $26$ standard benchmarks with different properties that cover most of the optimisation problems' nature, three traditional engineering problems, and one real complex engineering problem from the state-of-the-art literature. The algorithm performs well in finding global minima for nonconvex and multimodal functions, especially with high dimensional problems and it was found very competitive in comparison with the recent algorithmic proposals. Moreover, the algorithm outperforms and scales better than recent algorithms when it is benchmarked under a limited number of iterations for the composite CEC2017 problems. The design of this algorithm is kept simple so it can be easily coupled or hybridised with other search paradigms. The code of the algorithm is provided in C++14, Python3.7, and Octave (Matlab).
\end{abstract}

\section{Introduction}
\label{sec:introduction}
Over the past few decades, global optimisation techniques for solving combinatorial problems have flourished. The ever-increasing complexity of engineering applications means that variable sets are growing larger, and the subsequent landscapes that need to be explored by these optimisation problems are becoming increasingly complicated.

Many metaheuristic algorithms were developed to solve optimisation problems by mimicking biological and physical analogies \cite{REF:19}, including algorithms such as Genetic Algorithms (GAs) \cite{REF:12}, Particle Swarm Optimisation (PSO) \cite{REF:11} and the recent generalized version (GEPSO) \cite{REF:37}, and Harmony Search (HS) \cite{REF:13} and its latest modifications such as the \cite{REF:43, REF:1a, REF:44}, as well as the recent Whales Optimisation Algorithm (WOA) \cite{REF:6}, and Pathfinder Algorithm (PFA) \cite{REF:7}, to mention but a few. The candidate problems usually range from continuous differentiable problems to discrete, noisy, and even loosely defined objectives, such as in engineering applications.

Sequential sampling was used intensively in the second half of the last century for estimating the distribution of unknown functions, and the influence of such methods is still highly echoed in current engineering metamodels \cite{REF:24}. As the purpose of this paper is not to review all the optimisation algorithms that employed sequential sampling techniques, we advise the readers to refer to Homem-de-Mello and Bayraksan (2014) \cite{REF:47} where they can find a comprehensive review for algorithms used the Monte Carlo sampling approach for solving optimisation problems. Holistically, metaheuristics can be seen as performance-driven sequential sampling processes (Markovian processes) \cite{REF:17}, where the efficiency of such algorithms depends on the random walks that are used to explore and exploit the provided landscapes. Many random walks have been proposed in the literature, such as Brownian motion (obeys Gaussian distribution) and the L\'evy flights (obeys L\'evy distribution) for handling stochastic optimisation problems (see Yang et al. \cite{REF:17, REF:15} for comparison).

Estimation of distribution algorithms (EDAs) is another rank of search methods that tries to obtain better solutions from the search domain by estimating the probability density function (PDF) over the whole landscape from the so-far sampled solutions (see \cite{REF:28} for a detailed review). One relatively recent method that has been used by many engineering problems was proposed by Raphael and Smith (2003) \cite{REF:26}. In their method, they subdivide the search domain into a fixed number of sub-domains and depending on the fitness of the collected solutions the PDF of the domain evolves over time and the solutions are eventually focused around the most promising regions.

Generally, most current memetic algorithms are population-based and heavily dependent--as an initial step--on crude random sampling. This initial population step is crucial for the overall performance of the algorithms, affects all subsequent search phases \cite{REF:19}, is independent of the objective function and is considered an unbiased step of metaheuristics. Said otherwise, the quality of the chosen solutions in the initial population depends only on the design-of-experiment (DOE) methods used (space-filling techniques). Such methods assure equal accessibility to each part of the landscape. A plethora of publications dealt with the influence of the initial sampling and the population size on the overall performance of metaheuristics (such as in Polkov and Bujok (2018) \cite{REF:21}). For the here proposed algorithm, we use at every iteration DOE methods for sampling the new generation, as this guarantees simplicity and minimum structural bias of the algorithm.

The algorithm that we propose in this paper builds on the \textit{gbest} topology that was originally used in the PSO algorithm by Eberhart and Kennedy (1995) \cite{REF:14}. This simple and conservative topology prevents the algorithm from exploiting too much information from the candidates, which could drag the algorithm into a structural bias (see Kononova et al. \cite{REF:16} for definition) or even a premature convergence in late stages. In this work, we demonstrate how this simple algorithm can perform surprisingly well for complex optimisation problems, sometimes outperforming the existing state-of-the-art metaheuristics.

We further postulate our motivation behind this algorithmic proposal where we will use it as a heuristic for solving physical packing problems. The chosen topology of this algorithm was tailored to analogise the exact mathematical solvers used in packing regular items (boxes) in containers. Those solvers use the domain-reduction techniques of the available packing domain as we add more items to the container; some areas become more occupied than others and exploring the search domain of such areas is not computationally efficient. However, this algorithm will be used for packing physical items that are highly irregular and nonconvex in containers. Such complex geometries make the use of traditional methods impractical. For a comprehensive review about the techniques of domain-reduction in mathematics and optimisation problems we refer the reader to \cite{REF:41, REF:39} and some examples such as \cite{REF:38} and for the packing optimisation problems literature see \cite{REF:40, REF:42}.

Our algorithm uses only the best solution from the previous generation. This best solution will be the center of the tightened search domain, which is referred to as the current promising region. The size of the current promising region will be a function of time--the more optimisation steps that have already been completed, the smaller the promising region will be. The boundaries of the prominent region are only updated if a better solution than any previously obtained solution has been found. We sample inside and outside this region using a standard DOE method; in this paper we use Monte-Carlo sampling, though any other DOE method could also be used such as the Latin Hypercube Sampling (LHS) method (see \cite{REF:50, REF:51, REF:52}).

The chosen topology in this paper could be good for global exploration, as it maintains the population diversity, i.e., avoids premature convergence. Conversely, it could be problematic for exploiting local solutions where this process normally depends on the intensified information flow among the members of the population. Thus, this algorithm might be best hybridised in sequential or parallel schemes with other well-known algorithms. For more about topologies in metaheuristics, we refer the readers to \cite{REF:45} as they comprehensively review the common algorithmic topologies used in the literature. However, in this paper, we use the proposed algorithm as a standalone metaheuristic, with its simple topology, when benchmarking the algorithm against a range of problems and engineering applications.

The renowned No Free Lunch (NFL) theorems by Wolpert and Macready (1997) \cite{REF:2} have proven the non-uniformity of the superiority of algorithms for all ranks of optimisation problems. In this spirit, we propose that our algorithm be used as template, which can be combined with various DOE methods and other optimisation algorithms.

This work has been organized such that Section \ref{sec:algorithm} explains the analogy, the topology and the mathematical formulation of this algorithm. Section \ref{sec:parameters} provides a simple probabilistic analysis of the parameter settings, including a numerical illustration. In Section \ref{sec:benchmarks}, we test and compare the performance of this simplified approach with the state-of-the-art algorithms. Finally, our conclusions are provided in Section \ref{sec:conclusions}, some future applications for this approach are given in Section \ref{sec:future} and links to the source code of the algorithm are provided in Section \ref{sec:reproducibility}.

% ********* algorithm / topology *********
\section{Pareto-like sequential sampling (PSS)}
\label{sec:algorithm}
\par
\subsection{Analogy}
As mentioned earlier, the main analogy of this model is to mimic the Pareto principle, which implies that a \textit{"virtual few"} of a population account for the majority of the effects \cite{REF:3,REF:4}. The Pareto principal is popularly referred to as the 80/20 rule. This principle has been applied to a wide range of human relations and can be observed when a small percentage of contributors make the bulk impact on human-related matters. Historically, this 80/20 principle was first observed and described by the Italian professor Vilfredo Pareto from the University of Lausanne (UNIL), Switzerland (Pareto 1897) \cite{REF:4}.

We here employ the Pareto principle by more densely sampling the solutions' design variables from a tightened search domain, i.e., the current prominent region, while sampling the remainder of the solutions from the overall search domain.

In order to better understand the collective behavior considered in our analogy, one out of many possible examples, is by letting us imagine a selected group of people were appointed to discuss a certain social problem in a local community to find a proper solution. First, the problem will be communicated to each one of them individually and each one of them will have an initial opinion about what the solution could look like. However, these appointees, as representatives for the society, must sit together and meet regularly in order to reach a consensus on the optimal solution. As the appointees regularly meet, every one of them must express her/his own opinion, perspective, and exchange ideas of different solutions. Every meeting the \textit{majority} of these individuals will group behind certain revealing ideas that could be better to settle down the problem and try even to push it further by combining different solutions or sub-solutions from other individuals to reach an optimal solution. While the rest will always try to find better and optimal ideas that could convince the majority to shift to their sides. In this particular example, the appointees are the population of the solutions as they start randomly at one point (before the meetings and when they got informed) and then they refine their ideas every time they meet (number of iterations).

\subsection{Topology}
\begin{figure*}[htbp]
\centerline{\includegraphics[width = 1.3\linewidth]{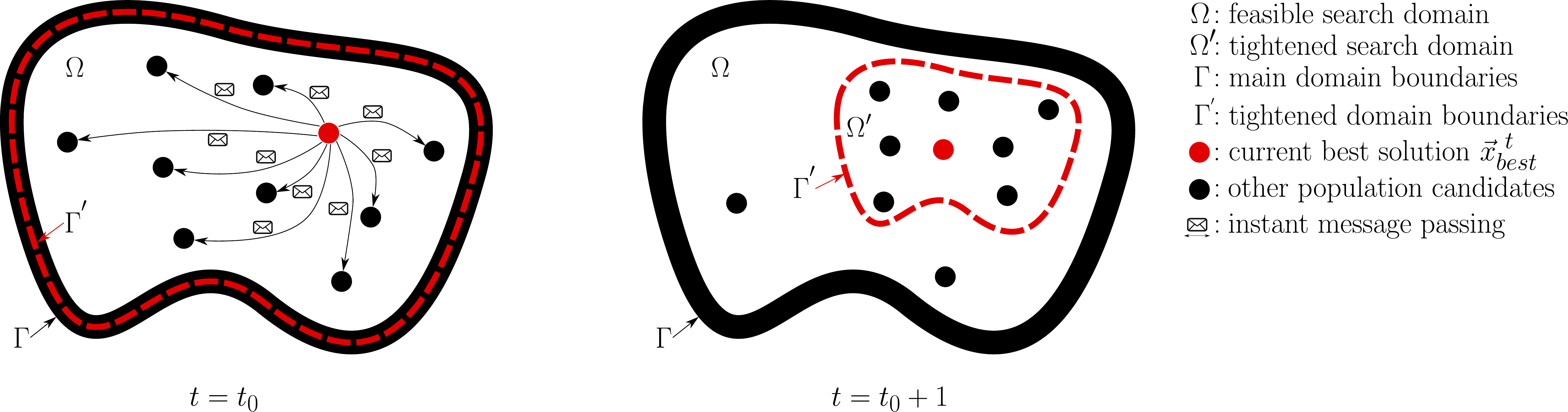}}
\caption{The topology of Pareto-like Sequential Sampling (PSS) algorithm.}
\label{FIG:1}
\end{figure*}

In this paper, we achieve the proposed analogy using mere performance-driven sequential sampling (see Fig. \ref{FIG:1}) to create better generations every iteration. We do this in a way that is compatible with any sampling method, i.e., classical DOE methods, such as the traditional Monte Carlo (MC) and Latin Hypercube Sampling (LHS), though in this paper we only use the MC sampling approach. Unlike traditional Evolutionary Algorithms (EAs), this algorithm does not use operations such as mutation or crossover to obtain better generations.

The process builds on the assumption that the current best solution vector, $\vec{x}_{best}$, lays in a sub-domain where the majority of the \textit{good} features exist (components of $\vec{x}$). We refer to this sub-domain as the prominent region $\Omega^{'}$. Hence, we sample the majority of the features, approximately $\alpha \times 100\%$, from the surrounding neighborhood. This topology of keeping track of only $\vec{x}_{best}$ each iteration is similar to the traditional \textit{gbest} topology proposed by Eberhart and Kennedy \cite{REF:14}.

The location and the size of the prominent domain changes every time the $\vec{x}_{best}$ changes. The size of the prominent domain (bandwidth), as will be described later, depends on the time remaining for the algorithm to make improvements and the predetermined size of the analogy $(1-\alpha)$. The size of the prominent region is set to a maximum of $(1-\alpha)\%$ of $\Omega$, centered about $\vec{x}_{best}$, where $\alpha$ is the acceptance probability that allows us to sample the components of the solutions from the prominent neighbourhood. The closer the algorithm gets to the maximum number of iterations, i.e., the shorter the time left for the algorithm to make improvements, the smaller the size of the prominent region.

This algorithm \textit{explores} and \textit{exploits} the domain by giving special emphasis to the prominent region around the current best solution. At each iteration, the algorithm performs the following steps:
\begin{itemize}
    \item Determine the new best fitness $\vec{x}_{~best}^{~i}$, 
    \item Update the prominent region $\vec{\eta}^{~i}$ if $~\vec{x}_{~best}^{~i}$ is better than $\vec{x}_{~best}^{~i-1}$,
    \item Sample from the search domain using classical DOE methods; sample more densely from the prominent regions $\vec{\eta}^{~i}$ and less densely from the overall search domain.
\end{itemize}

\subsection{Mathematical Model}
\label{subsec:mathmodel}
Let the real-valued optimization problem $f(\vec{x})$ be defined on the search domain $\Omega \in {\rm I\!R}^{n}$ and for all sub-domains $^{i}\Omega^{'} \subset \Omega$. The corresponding boundaries are $\vec{\Gamma}$ and $^{i}\vec{\Gamma^{'}}$ for the main domain and the sub-domain, respectively.
\begin{equation}
    \text{optimize}~f(\vec{x}), ~{\vec{x}} \in \Omega, \Omega \subset {\rm I\!R}^{n},~ \text{where}
\label{eqn:1}
\end{equation}
$$\vec{x} = \{x_1, \dots, x_n\},~ \vec{\Gamma}^{-}\leq \vec{x} \leq \vec{\Gamma}^{+}.$$
More specifically, the lower and the upper bounds for the main domain and the sub-domain can be written as $[\vec{\Gamma}^{-}, \vec{\Gamma}^{+}$] and $[ ^{i}\vec{\Gamma}^{'-}, ^{i}\vec{\Gamma}^{'+}]$, respectively. Additionally, let there be a solution vector $\vec{x}_{opt}$ that reclines in $\Omega^{'}$ and minimizes $f(.)$, where $f(\vec{x}_{opt}) \leq f(\vec{x})~ \forall~ \vec{x} \in \Omega$ is a global solution. For clarity, the following notations will be used in this work: $i$ is the current time step (iteration), $k$ is the index of the solution vector in the population matrix and $j$ indicates the $j^{th}$ design variable (feature $x_j$) of any solution vector $\vec{x}$.

The process of random feature selection depends either on the algorithm sampling a decision variable $x_{j}^{i} \in \vec{x}_j$ from the prominent domain $\Omega^{'}$ or improvising and sampling randomly from the entire region $\Omega$. This decision is dependent on the acceptance probability $\alpha$.

Like any other metaheuristic algorithm, the first step is the initial population sampling. For a continuous landscape, the population can be composed by:
\begin{equation}
    _{k}\vec{x} = \vec{\Gamma}^{-} +~ _{k}\vec{u} \odot(\vec{\Gamma}^{+} - \vec{\Gamma}^{-}).
\label{eqn:2}
\end{equation}
In this equation, $_{k}\vec{u}$ is a vector of random coefficients sampled using one of the classical DOE methods. The $\odot$ operator indicates an element-wise multiplication of vectors.

The matrix $\mathbf{u}^{i}$ is sampled at each iteration $i$ with the chosen DOE method. It contains the sampled random coefficients; its size is $\beta \times n$ where $\beta$ is the size of the population and $n$ is the number of dimensions of the problem.
\begin{equation}
[\mathbf{u}^{i}] =
\begin{bmatrix}
_{1}u_{1}^{i} & _{1}u_{2}^{i} & \dots & _{1}u_{n}^{i} \\ 
_{2}u_{1}^{i} & _{2}u_{2}^{i} & \dots & _{2}u_{n}^{i} \\ 
\vdots & \vdots & \ddots & \vdots\\ 
_{\beta}u_{1}^{i} & _{\beta}u_{2}^{i} & \dots & _{\beta}u_{n}^{i} \\ 
\end{bmatrix}.
\label{eqn:3}
\end{equation}

After evaluating the initial population and determining the best solution $\vec{x}_{~best}^{~i}$, the algorithm estimates the domain-tightening coefficients to update the current prominent upper and lower bounds for each design feature as per the following expressions:
\begin{equation}
    ^{i}\vec{\Gamma}^{~'+} = \vec{x}_{~best}^{~i} + \vec{\eta}^{~i},~ ^{i}\vec{\Gamma}^{~'+} \leq~ \vec{\Gamma}^{+},
\label{eqn:4}
\end{equation}
\begin{equation}
    ^{i}\vec{\Gamma}^{~'-} = \vec{x}_{~best}^{~i} - \vec{\eta}^{~i},~ ^{i}\vec{\Gamma}^{~'-} \geq~ \vec{\Gamma}^{-},
\label{eqn:5}
\end{equation}
\begin{equation}
    \vec{\eta}^{~i} = \frac{(1-\alpha)\Big( 1 - \frac{i}{\gamma}\Big)}{2}~ (\vec{\Gamma}^{~+} -~ \vec{\Gamma}^{~-}).
\label{eqn:6}
\end{equation}

The bandwidth ($\vec{\eta}^{~i}$) is updated every time $\vec{x}^{~i}_{~best}$ is updated (self-adaptive mechanism). As can be seen in (\ref{eqn:6}), we applied time-dependent bandwidth tightening using the term $\big(1-\frac{i}{\gamma}\big)$. This controls the greediness of the algorithm with time: as time passes, the algorithm becomes greedier. The $(1-\alpha)$ term in (\ref{eqn:6}) is used to maintain the $\alpha / (1-\alpha)$ ratio (resembling the 80/20 analogy). To control the number of features (design variables) that will be sampled per solution from the prominent domain $\Omega^{'}$, the acceptance probability $\alpha$ has been used.

If the algorithm decides to draw from the prominent domain $\Omega^{'}$, the set of equations (\ref{eqn:4}) to (\ref{eqn:6}) together with (\ref{eqn:7}) are used to generate a feature as follows:
\begin{equation}
    _{k}{x}_j^i =~ ^{i}_{j}{\Gamma}^{'-} +~ _{k}{u}_j^i ~(^{i}_{j}{\Gamma}^{'+} -~ ^{i}_{j}{\Gamma}^{'-}).
\label{eqn:7}
\end{equation}
If the algorithm decides to draw a feature from the overall domain $\Omega$, then it will instead be evaluated using:
\begin{equation}
    _{k}{x}_j^i =~ _{j}{\Gamma}^{-} +~ _{k}{u}_j^i ~(_{j}{\Gamma}^{+} -~ _{j}{\Gamma}^{-}).
\label{eqn:8}
\end{equation}
This recursive process runs until the stopping condition is satisfied, which we have set to the maximum number of evaluations $\gamma$. Algorithm \ref{alg:main} explains the steps of the algorithm.

\begin{algorithm}
\SetAlgoLined
\KwResult{Minimize the objective function $f(\vec{x})$}
 read parameter-settings\;
 pre-allocate and initialize memory containers\;
 \Begin{
 initialize a population sample $\in \vec{\Gamma}$\;
 evaluate the initial population\; %and determine $\vec{x}_{best}^{~0}$
    \For{$i = 1~$\KwTo$\gamma$}{
        \tcp{best solution in iteration $i$ is $\vec{x}^{~i}_{best}$}
        find the current best $\vec{x}^{~i}_{best}$\;
        generate $[\mathbf{u}^i]_{\beta \times n} \sim U(0,1)$\;
        \For{$k=1~$\KwTo$\beta$}{
            \If{$ f(\vec{x}^{~i}_{best}) < f(\vec{x}^{~i-1}_{best})$}{
                update $\vec{\eta}^{~i}$\ and $^{i}\vec{\Gamma}^{'}$\;
            }
            \ForEach{$x_j \in \vec{x}$}{
                \uIf{$ rand \sim U(0,1)\leq \alpha$}{
                    \tcp{sample from the prominent domain $\Omega^{'}$}
                    choose random $_{k}x^{i}_{j} \in \Gamma^{'}_{j}$;
                }
                \Else{
                    \tcp{sample from the overall domain $\Omega$}
                    choose random $_{k}x^{i}_{j} \in \Gamma_{j}$;
                }
            }
            $[population]_{\beta \times n} \leftarrow~ _{k}\vec{x}$\;
            evaluate $f(_{k}\vec{x})$\;
        }
    }
}
    \Return{$\vec{x}_{best}$}\;
 \caption{Pareto-like Sequential Sampling (PSS)}
 \label{alg:main}
\end{algorithm}

% ********* Parameter tuning *********
\section{Parameters Setting}
\label{sec:parameters}
\par
\begin{figure}[htbp]
\centerline{\includegraphics[width = 0.4\linewidth]{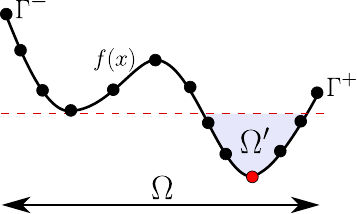}}
\caption{A virtual 1D optimization objective $f(\vec{x})$ -- discretized space.}
\label{FIG:5}
\end{figure}

The No-Free-Lunch (NFL) theorem directly promotes that no algorithm has a global parameters set that ensures an overall good performance for all the possible optimisation problems \cite{REF:29}. In line with this, we provide a comprehensive theoretical and visual explanation of the parameter-setting (tuning) for the PSS algorithm. A survey about traditional tuning techniques can be found in \cite{REF:49, REF:23}. The PSS algorithm is based on three input parameters that need to be set by the user. Namely, the population size ($\beta$), the maximum number of iterations ($\gamma$) and the acceptance probability ($\alpha$). The taxonomy and terminology used in this section are taken from Eiben et al. \cite{REF:23}. The derived formulae here follow probabilistic rules that can be found in standard textbooks about the theory of probability (for example see \cite{REF:55}). In this section, we provide an analytical investigation of the interaction between the three parameters and illustrate their effects on finding the global optimum by means of a numerical example. Apart from this, the self-adaptive mechanism for the prominent domain will not be discussed here (see Section \ref{subsec:mathmodel}).

Let $f(x)$ be a 1D discrete optimization problem with $N$ finite possible solutions that we need to minimize. Fig. \ref{FIG:5} shows the distribution of $x$ along the domain $\Omega$. Now, let $\Omega^{'}$ be the \textit{true}\footnote{This is an innate constitutive property of the distribution (the change) of $f(\vec{x})$ vs. any independent $x_j \in \vec{x}$.} prominent domain around the optimum solution $x_{opt}$. The red-dashed line, shown in Fig. \ref{FIG:5}, represents the lowest valley of all local minima, wherein all the points below this line are better than any local optima above it. If any solution below this line is chosen, it will always lead to the global solution in the coming generations by using the \textit{gbest} topology. Assuming that $n^{'}$ is the number of all possible points below the line $\in \Omega^{'}$ and $N$ is the number of all possible solutions $\in \Omega$, the probability of finding a solution inside the true prominent domain in the first step (initial population sampling) is:
\begin{equation}
    P(x \in \Omega^{'}) = \dfrac{n^{'}}{N}.
\label{eqn:9}
\end{equation}
We can expand (\ref{eqn:9}) to determine the probability of drawing at least one solution that lays in $\Omega^{'}$ (event $A_0$) from $\beta$ drawn solutions:
\begin{equation}
    P(A_0) = 1 - \bigg(1 - \dfrac{n^{'}}{N}\bigg)^{\beta}.
\label{eqn:10}
\end{equation}
Equation (\ref{eqn:10}) explains the exponential relationship between the landscape size and the initial population size. The probability of event $A_0$ is a problem-specific property that differs per objective and independent design variable. The complementary event to $A_0$ will be called $B_0$, wherein $P(B_0) = 1 - P(A_0)$ and $P(B_0)$ means none of the drawn $\beta$ solutions initially is located in the \textit{true} prominent region.

After the population initialization, the algorithm performance phase begins. In this phase, the algorithm tests the acceptance-rejection paradigm against $\alpha$ for each solution component ($\forall x^i_j \in \vec{x}^i$). If a random number $r^i\sim U(0,1) \leq \alpha$, the algorithm will sample the design component from the prominent region, else it will be sampled from the overall region. The probability of exploring new solutions $\in \Omega^{'}$ from the domain is defined in Equation (\ref{eqn:11}). This equation is independent of the time step, $i$, and it is always constant if the population size $\beta$ is constant for a specific optimization problem and an independent design variable.
\begin{equation}
    P^i\big(x\in \Omega^{'}|r^{i} > \alpha\big) = \left[  1 - \bigg(1 - \dfrac{n^{'}}{N}\bigg)^{\beta} \right] (1-\alpha).
\label{eqn:11}
\end{equation}

Now if the algorithm were to intensificate (success against $\alpha$), the probability of getting a better solution within the current fictitious prominent region is shown in the following:
\begin{equation}
    P^i\big(x^{i+1}_{best} < x^{i}_{best}|r^{i} \leq \alpha\big) = \left[  1 - \bigg(1 - \dfrac{n_{t}^{i}}{N_{t}^{i}}\bigg)^{\beta} \right] \alpha.
\label{eqn:12}
\end{equation}

Notice that Equation (\ref{eqn:12}) is time-dependent. The terms $n_t^{i}$ and $N_t^{i}$ by definition depend on the sub-region (the current \textit{fictitious} prominent region), and they are a function of time (notice the $i$ superscription). This can be imagined as if the red-dotted line corresponds to the current best solution and the algorithm tries to find a better sub-optimal. Put differently, it implies--only in this algorithm--that the intensification property is time-dependent and it resembles a perturbation problem.

In real life, the algorithm is not aware of the \textit{true} prominent domain. Indeed, the only parameters known to the algorithm are the domain and the objective function (as an evaluator). We assume that the best-known solution so far lays in the prominent domain and hope that the algorithm's exploration will find a better one in the upcoming generations (if any).

If we have an optimization problem with $n$ dimensions where $n > 1$, we always assume independent design variables. For expanding the reflections made above for $n$ design variables, the probability of obtaining an optimal solution can be calculated using the multiplication rule, wherein every design variable is treated as an independent event. This simple probabilistic analysis shows that the algorithm is dependent on the distribution of the objective function, though because the algorithmic parameters are interactive, they must be carefully set. From this analysis, we can see that with a bigger $\beta$, we boost the probability of obtaining solutions that belong to the true prominent region, i.e., global optimum. With a higher $\alpha$, the algorithm puts most of its effort into intensification, weakening the diversification. For multimodal and nonconvex problems, $\alpha$ should be smaller than in convex and unimodal problems. Using more iterations will always increase the probability of exploring better solutions. Finally, it is worth mentioning that we used only MC sampling in the presented analysis. Unless mentioned otherwise, we used a population of $30$ and an acceptance probability of $0.95$.

\subsection{Numerical Illustration}
\begin{figure*}[htbp]
\centerline{\includegraphics[width = 1.5\linewidth]{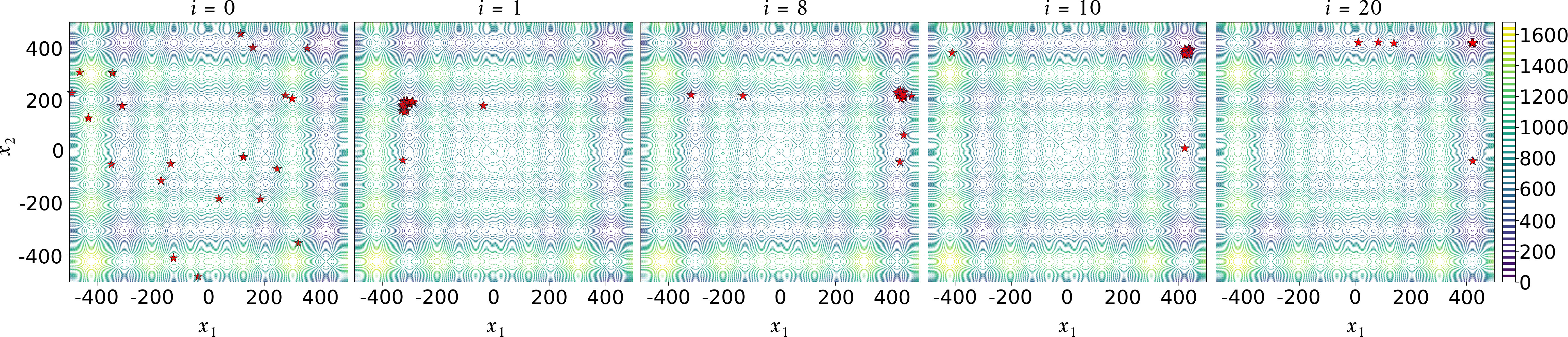}}
\caption{The population history for Case (a) -- 2D Schwefel function.}
\label{FIG:3}
\end{figure*}

\begin{figure}[htbp]
\centerline{\includegraphics[width = 0.60\linewidth]{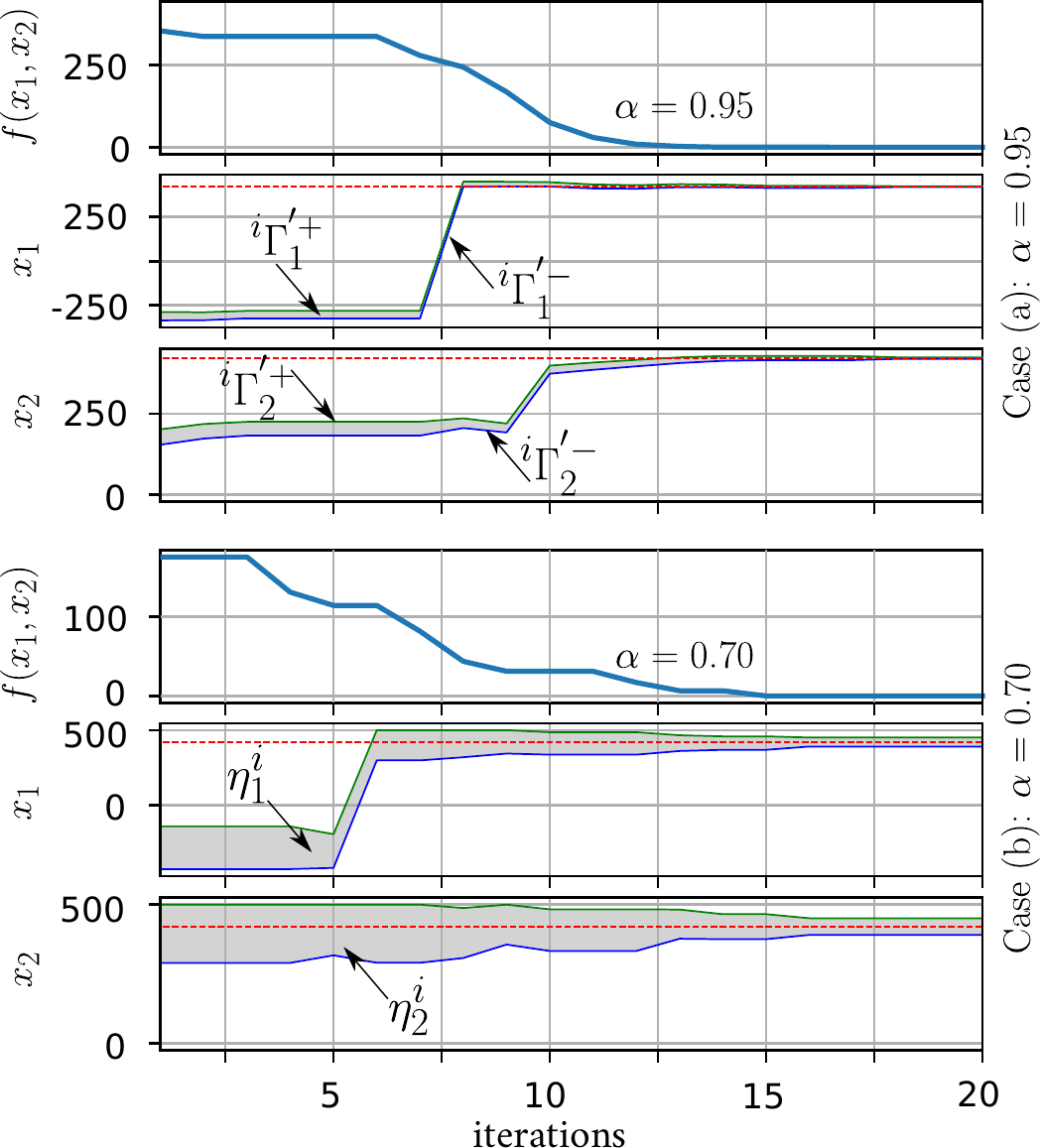}}
\caption{Domain tightening and convergence histories -- 2D Schwefel function.}
\label{FIG:4}
\end{figure}

The example shown in Fig. \ref{FIG:3} and \ref{FIG:4} indicates the solution of the standard 2D Schwefel problem. The optimum solution $f(\Vec{x}_{opt}) = 0.0$ is obtained when both $x_1$ and $x_2$ hit $420.9687$ (the red-dotted lines in Fig. \ref{FIG:4}). A population of $30$ members and a sum of $20$ iterations were used to solve this problem. For this problem, we set $\alpha = 0.95$ and $0.70$. Fig. \ref{FIG:4} reveals how the algorithm adaptively changed the prominent search area $\Omega^{'}$ (the gray fill-up) and dynamically tightened it to control the greediness of the algorithm with respect to the remaining time (for the two cases).

\begin{figure*}[htbp]
\centerline{\includegraphics[width = 1.3\linewidth]{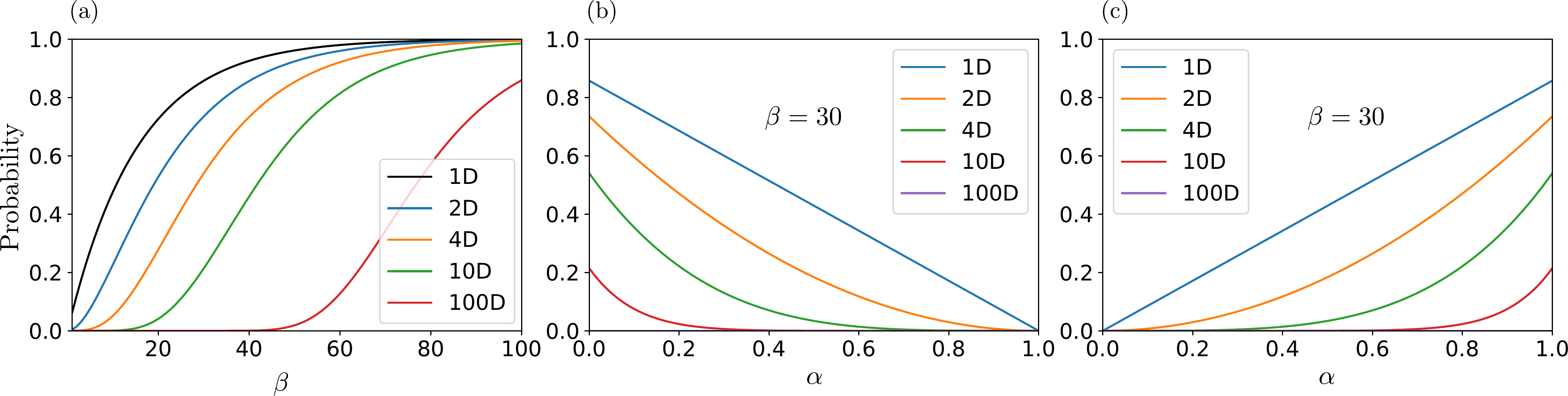}}
\caption{(a) The change of $P(A_0)$ with $\beta$; (b) the diversification probability vs. $\alpha$ in $\Omega$; (c) the intensification probability vs. $\alpha$ in $\Omega^{'}$ -- $n$D Schwefel function for $n=1$, $n=2$, $n=4$, $n=10$, and $n=100$.}
\label{FIG:4b}
\end{figure*}

According to our definition of the true prominent region for Schwefel function $\Omega^{'} \in [389.33, 452.16]~ \forall~ x_j$ (determined graphically), we say that the algorithm has globally converged if the algorithm successfully found both $x_1$ and $x_2 \in \Omega^{'}$. We analyzed $30$ consecutive runs for each case, and the converged results are expressed as: average $\pm$ standard deviation (success rate \%). For Case (a), the results were 0.197345 $\pm$ 0.835687 (83.33\%), and for Case (b), the results were 2.435370 $\pm$ 2.599124 (96.67\%). These results highlight the trade-off balance between the diversification and the intensification processes explained in Equations (\ref{eqn:11}) and (\ref{eqn:12}). Side by side, Fig. \ref{FIG:4b} (a) reveals the effect of choosing the population size ($\beta$) on the probability of finding better initial solutions (event $A_0$). Moreover, Fig. \ref{FIG:4b} (b) and \ref{FIG:4b} (c) explain the global diversification ($\forall x_i \in \Omega$) and the intensification inside the true prominent domain ($\forall x_i \in \Omega^{'}$) for different dimensions of Schwefel function, respectively.

% ********* Results and Discussions *********
\section{Benchmarks and Comparisons}
\label{sec:benchmarks}
\par
The presented PSS algorithm was used to solve a set of standard functions (Section \ref{subsec:stdbenchmarking}), the CEC2017 composite functions (Section \ref{subsec:compo}), and another selected set of engineering problems (Section \ref{subsec:engbenchmarking}). We benchmarked and compared the performance of the PSS algorithm with the obtained results from other state-of-the-art algorithms, the results of which were mostly retrieved directly from the source papers and not re-simulated herein. To obtain a meaningful comparison, we used the same number of evaluations per problem as were used by the authors of the other algorithms. This benchmarking approach follows the recommendations by Arcuri and Fraser (2013) \cite{REF:29} and Liao et al (2015)\cite{REF:30}, which were recently reemphasized by Ser et al. (2019) \cite{REF:19} and used by Piotrowski and Napiorkowski (2018) \cite{REF:20}.

For the standard benchmarks, we chose state-of-the-art algorithms to compare with and avoided using the classical and outperformed algorithms following the suggestions by \cite{REF:19,REF:31}. The first algorithm we used in the benchmarking was the Whale Optimization Algorithm (WOA) as this is one of the recent and heavily cited algorithms in the literature, which outperformed many other recent algorithmic proposals (see Mirjalili and Lewis (2016) \cite{REF:6}). The second one was the Pathfinder Algorithm (PFA), which has also been published recently and shown to outperform many other algorithms (see Yapici and Cetinkaya (2019) \cite{REF:7}).

For engineering benchmarks, we compared the obtained solution by directly adopting the best results obtained from different algorithms in the literature. We also solved a recent engineering case study with high dimensions that was solved by using the Modified Parameter-Setting-Free Harmony Search (MPSFHS) algorithm by Shaqfa and Orb\'an (2019) \cite{REF:1a} and we used the same algorithm to benchmark the new algorithm with regard to its scalability.

\subsection{Standard Benchmarks}
\label{subsec:stdbenchmarking}
Here, we compare our proposed algorithm with the WOA \cite{REF:6} and the PFA \cite{REF:7} using the standard benchmarks explained in Table \ref{tab:benchmarks}. We carefully chose the benchmarks in this paper to cover a wide range of problems.

In Table \ref{tab:WOA}, we compared the WOA with our proposed PSS method. The problems were solved by assuming $\alpha = 0.95$. For the population size and the number of iterations, we chose the same values as in \cite{REF:6}, i.e., a population size of $30$ and a total of $500$ iterations. In general, $30$ dimensions were used for all the proposed problems ($n=30$) except for $f_{13}$ and $f_{15}$ (see Table \ref{tab:benchmarks}).

The reported results in Table \ref{tab:WOA} express the average $\pm$ the standard deviation for $25$ consecutive runs per problem (as in \cite{REF:6}). The boldfaced results indicate--for each problem--the algorithm that performed the best in the comparison. From Table \ref{tab:WOA}, it can be seen that in easy unimodal and convex functions, the best algorithm was WOA. This behaviour was expected due to the weak intensification in the PSS algorithm that the crude random walk usually reveals. For harder problems, such as nonconvex and/or multimodal ones, though, the PSS algorithm was better able to allocate global optima--this finding holds notably for the Schwefel function $f_{11}$.

In Table \ref{tab:PFA}, we compare the PSS algorithm to the recently proposed PFA algorithm \cite{REF:7}. As in their work, the problems in this table were simulated using a population size of $30$ and $1000$ iterations. For the PSS algorithm, $\alpha$ was set to $0.95$ as per the previous comparison. However, the number of dimensions used for $f_{8}$ was $20$, while for $f_{7}$ $n=6$, and for $f_{14}$ and $f_{16}$, $n$ was set to $2$ and $3$, respectively, as shown in Table \ref{tab:benchmarks}. The stated results are for $30$ consecutive runs per function (as in \cite{REF:7}).

As can be seen in Table \ref{tab:PFA}, $f_2$ and $f_5$ behaved the same as in Table \ref{tab:WOA} (the pronounced weak intensification for the proposed approach). On the other hand, the results of $f_{8}$ and $f_{11}$, see Table \ref{tab:PFA}, were outperformed by the Pareto-like sampling. The results of $f_{7}$, $f_{14}$, and $f_{16}$ are almost identical.

\subsection{Testing Scalability}
\label{subsec:scalability}
With an increased dimensionality, the problem complexity increases, and accordingly, the search domain expands exponentially. To test the scalability of the current algorithm, we compared the behaviour of the proposed approach with the PFA algorithm for $f_8$ and $f_{11}$ functions. The simulations were run $30$ times for each, with $10$, $50$, and $100$ dimensions as illustrated in Table \ref{tab:PFAScale}. As the results suggest, the proposed PSS algorithm outperformed the PFA and registered fewer deteriorations in performance with the increase of dimensions.

As this test implies that increased dimensionality requires more computational capacity, we ran simulations with the recently modified parameter-setting-free harmony search algorithm (MPSFHS) and the proposed PSS approach to see how increasing the iterations could enhance the results. To do this, $f_{8}$ and $f_{11}$ were solved for $n=100$ and a total of $300,000$ evaluations ($10,000$ iterations with a population size of $30$ for PSS). Table \ref{tab:MPSFHS} illustrates the results of $30$ consecutive runs for each algorithm. The PSS approach outperformed the MPSFHS algorithm, and considerable enhancements were seen in the end results.

\subsection{Composite Benchmarks}
\label{subsec:compo}
In this section, we used the CEC2017 by Wu et al. (2016) \cite{REF:22} and reviewed by Molina et al. (2018) \cite{REF:48}, competition's composite functions to benchmark the behaviour of the proposed algorithm (PSS). The used functions are briefly described in Table \ref{tab:composite_benchmarks}. We first tested our algorithm against these functions under extreme cases where only a few iterations are allowed. In Table \ref{tab:PSS_COMPO}, we ran $30$ consecutive tests with $n = 2$, $\alpha = 0.95$, $\beta = 30$, and with allowing a total of $10$ iterations per run for the composite functions $f_{c,1} \to f_{c,8}$. We conducted a comparison under the same conditions for the PSS, PFA, WOA, and PSO algorithms. The MPSFHS was tuned to $HMS = 30$ $m=2$, $HMCR_i = 0.75$, $PAR_i = 0.15$, $HMCR_{max} = 0.99$,  $PAR_{min} = 0.05$, and $300$ iterations to have an equivalent number of evaluations as the PSS, PFA, WOA, and PSO (for more about the used parameter-settings refer to \cite{REF:1a}).

The results shown in Table \ref{tab:PSS_COMPO} suggest that the PSS algorithm can still be globally convergent and allocate global minima in most of the runs and it scales well with the available time (iterations). Indeed it scored the best results for most functions, though, the gap between the PSS and the PFA algorithms was small. The best runs are revealed in Fig. \ref{FIG:9} showing two iteration-milestones $i = 5$ and $i = 10$. Notice that the reported results in Table \ref{tab:PSS_COMPO} are expressed in terms of the error values and computed as ($f_{c,i}(\vec{x_{i}}) - f_{c,i}(\vec{x_{opt}})$) and we revealed the minimum, maximum, mean, median, and the standard deviation of the error values as recommended by the CEC2017 report (refer to \cite{REF:22}).

In another extreme case we tested the functions $f_{c,1} \to f_{c,10}$ but this time with $n=100$ and by only employing $100$ iterations. In this experiment, we only compared the PSS with the PFA algorithm. The shown results in Table \ref{tab:PSS_COMPO2} suggest that the gap between the PSS and the PFA widens. Indeed, the PSS outperformed the PFA in all the functions except $f_{c,2}$ and $f_{c,7}$. These results signify the capability of the PSS algorithm to scale with high-dimensional problems while exploiting very limited computational resources. This is quite important for complex surrogate models that require a significant computational capacity for each functional evaluation \cite{REF:53, REF:54}.

We also conducted simulations similar to the ones presented by Yapici and
Cetinkaya (2019) \cite{REF:7} and again we compared our results with the adopted ones directly from their paper. The tests were conducted with different dimensions for each function; we used $n=10$, $n=30$, and $n=50$. $50$ runs with $1000$ iterations per each were used to run all the tests. $\alpha$ and $\beta$ were also set to be $0.95$ and $30$, respectively. The results of both algorithms, shown in Table \ref{tab:compbench}, are in proximity to each other for most functions. However, in their paper they compared the performance of the PFA with other algorithms and the best performance was registered by the Effective Butterfly Optimizer (EBO) \cite{REF:34} and the enhanced version of the Success-History based Adaptive Differential Evolution LSHADE-cnEpSin algorithm \cite{REF:35} and they outperformed, though the results are not shown here (refer to \cite{REF:7}), the PFA and the PSS as those algorithms (EBO and LSHADE-cnEpSin) were specifically designed to solve the CEC2017 problems.

\begin{figure*}[htbp]
\centerline{\includegraphics[width = 1.0\linewidth]{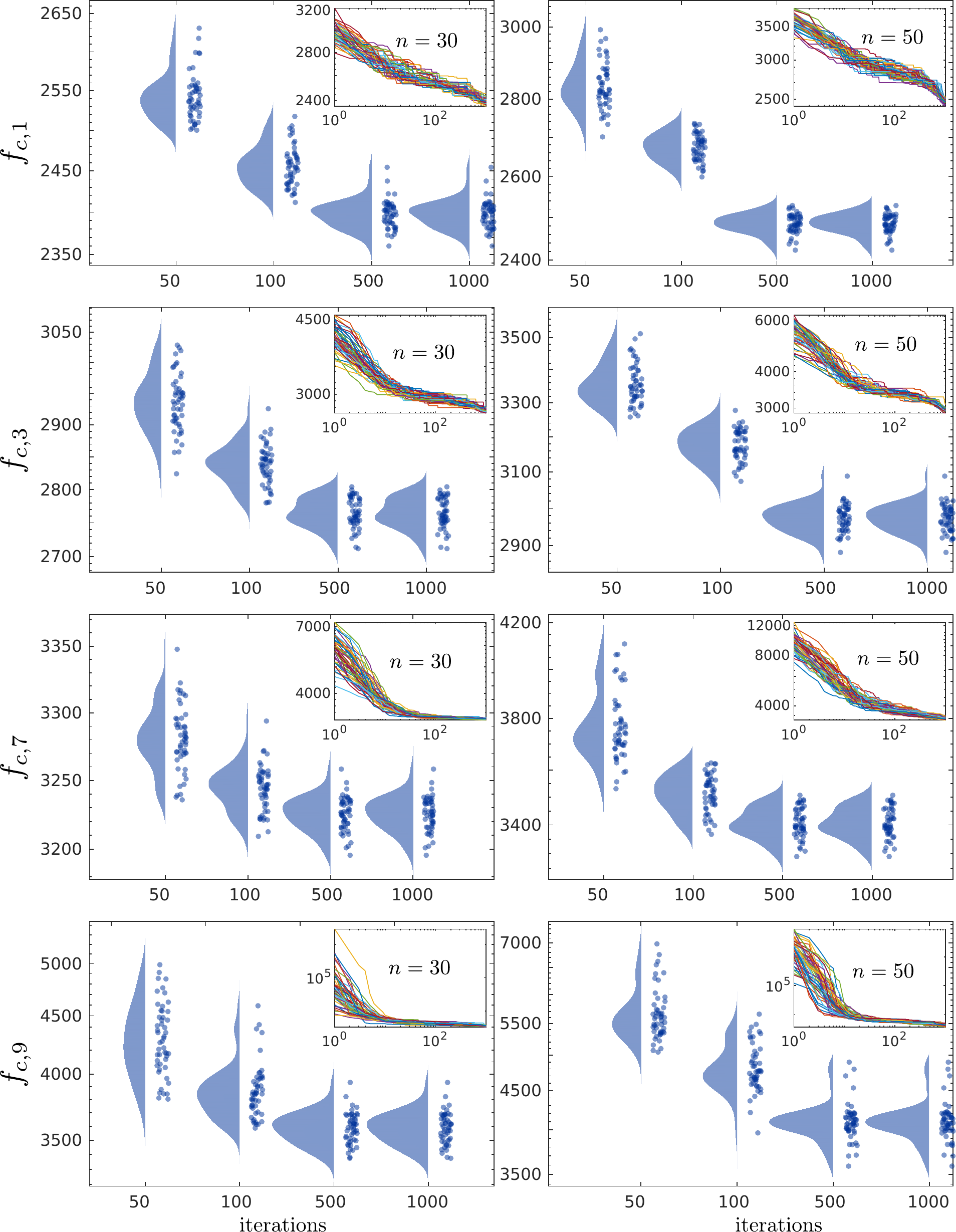}}
\caption{The statistical distribution of the obtained optima for a set of composite functions, with $50$ runs per each, computed at four iteration-milestones ($i=50$, $i=100$, $i=500$, and $i=1000$) with a semi-log scale (on the y-axis). (insets), reveal the actual log-log scale convergence curves for the all $50$ runs.}
\label{FIG:raincloud1}
\end{figure*}

The CEC2017 report indeed recommends to measure the performance of the algorithms at several iteration-milestones to monitor the convergence behavior. However, the reported error values in terms of the mean and standard deviation does not reveal the real distribution of the data. For this purpose, we used the raincloud plots by Allen et al. (2019) \cite{REF:36} to monitor the convergence at several iteration-milestones. In Fig. \ref{FIG:raincloud1}, we evinced a sample of the convergence data distribution at different milestones $i = 50$, $i = 100$, $i = 500$, and $i = 1000$ and for problem dimensions $n = 30$ and $n = 50$. As can be seen in the figure the results are widely-distributed at the beginning of the performance stage (notice data at $i = 50$); at later iterations relatively narrower distributions can be identified. It can be noticed as well that some distributions have two peaks and this could be related to the corresponding topology of the function as they contain many attractive local minima that could trick the algorithm into it (for example see $f_{c,7}$ in Fig. \ref{FIG:raincloud1} and the corresponding 2D surface in Fig. \ref{FIG:9}).

\begin{figure*}[htbp]
\centerline{\includegraphics[width = 1.0\linewidth]{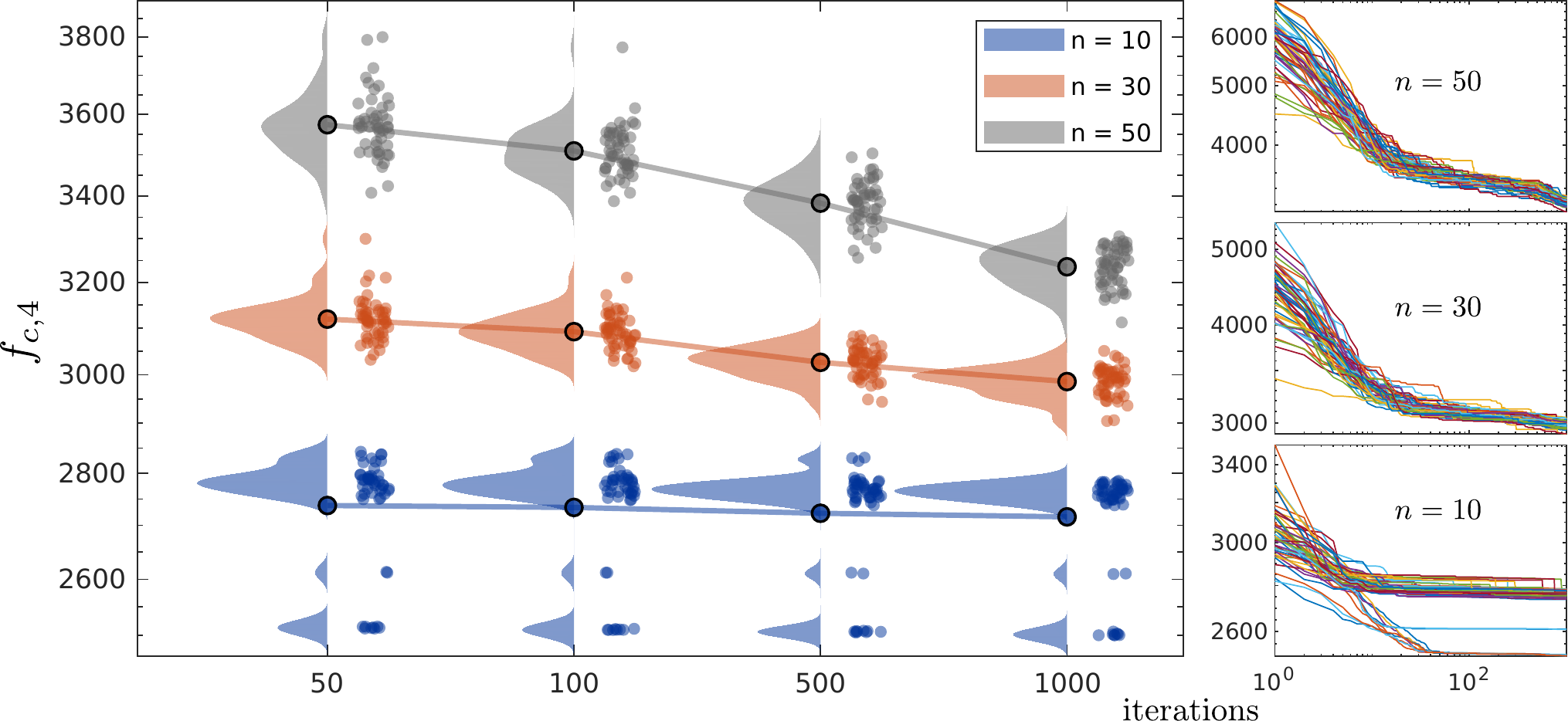}}
\caption{The statistical distribution of convergence data of the composite function $f_{c,4}$ for $n=10$, $n=30$, and $n=50$ at four different iteration-milestones ($i=50$, $i=100$, $i=500$, and $i=1000$). The right column shows the actual convergence curves in log-log scale corresponding to $n=10$, $n=30$, and $n=50$.}
\label{FIG:raincloud2}
\end{figure*}

The raincloud figures help to visualize how frequent the algorithm could be trapped in local minima and how often it converges globally. This is more clear in lower dimensions, for instance, when $n = 10$ as shown in Fig. \ref{FIG:raincloud2} three different distributions can be noticed as each one of them represents a possible local minima. The consistency of the algorithm can be seen as it generates more iterations it approaches global answers without being staggered at local minima. Moreover, the dimensional discrepancy can be clearly distinguished on the same figure for all the runs.

\begin{table*}[htbp]
\centering
\caption{Standard benchmark functions (refer to \protect\cite{REF:46} for details).}
\label{tab:benchmarks}
\begin{tabular}{ll}
\hline
\multicolumn{1}{c}{\textbf{Function}} & \multicolumn{1}{c}{\textbf{Expression}} \\\hline
\\\vspace{-5mm}\\
Sphere         &  $f_{1}(\vec{x}) = {\sum\limits_{i=1}^{n} x_i^{2}}
$\\
Sum squares    &  $f_{2}(\vec{x})=\sum\limits_{i=1}^{n}{ix_i^2}
$\\
Chung Reynolds &  $f_{3}(\vec{x}) = \Big({\sum\limits_{i=1}^{n} x_i^{2}}
\Big)^{2}$\\
Schwefel2.21   &  $f_{4}(\vec{x}) = \max\limits_{i=1,...,n}|x_i|$\\
Schwefel2.22   &  $f_{5}(\vec{x}) = \sum\limits_{i=1}^{n}|x_i|+\prod\limits_{i=1}^{n}|x_i|$ \\
Rosenbrock     &  $f_{6}(\vec{x})=\sum\limits_{i=1}^{n}[100 (x_{i+1} - x_i^2)^ 2 + (1 - x_i)^2]$\\
Trid 6         & $f_{7}(\vec{x}) = \sum\limits_{i=1}^{n}(x_i-1)^2 - \sum\limits_{i=2}^{n} x_i x_{i-1}$\\
 Zakharov       &  $f_{8}(\vec{x}) = \sum\limits_{i=1}^n x_i^{2}+(\sum\limits_{i=1}^n 0.5ix_i)^2 + (\sum\limits_{i=1}^n 0.5ix_i)^4$\\
Griewank       &  $f_{9}(\vec{x}) = 1 + \sum\limits_{i=1}^{n} \frac{x_i^{2}}{4000} - \prod\limits_{i=1}^{n}cos(\frac{x_i}{\sqrt{i}})$\\
Ackley         &  $f_{10}(\vec{x}) = -20exp\bigg(-0.2\sqrt{\frac{1}{n}\sum\limits_{i=1}^{n}x_i^2}\bigg)-exp\Big(\frac{1}{n}\sum\limits_{i=1}^{n}cos(2\pi x_i)\Big)+ 20 + exp(1)$\\
Schwefel       &  $f_{11}(\vec{x}) = 418.9829n -{\sum\limits_{i=1}^{n} x_i sin(\sqrt{|x_i|})}$\\
Shubert        & $f_{12}(\vec{x}) = \prod\limits_{i=1}^{n}{\left(\sum\limits_{j=1}^5{ cos\big((j+1)x_i+j\big)}\right)}$\\
%  & -$exp\Big(\frac{1}{n}\sum\limits_{i=1}^{n}cos(2\pi x_i)\Big)+ 20 + exp(1)$\\
Six-hump camel &  $f_{13}(x_1, x_2) = \Big( 4 - 2.1x_1^2 + \dfrac{x_1^4}{3}\Big)x_1^2 + x_{1} x_{2} + \Big(-4 + 4 x_2^2\Big)x_2^2$\\
Goldstein      & $f_{14}(x_1,x_2)=[1 + (x_1 + x_2 + 1)^2(19 - 14x_1 +3x_1^2- 14x_2 + 6x_1 x_2 + \dots$\\
 & $ 3x_2^2)][30 + (2x_1 - 3x_2)^2(18 - 32x_1 + 12x_1^2 + 48x_2 - 36x_1 x_2 + 27x_2^2)]$\\
De Jong 5      & $f_{15}(x_1, x_2) = \bigg(\dfrac{1}{500} + \sum\limits_{j=1}^{25}\Big( j + \sum\limits_{i = 1}^{2}(x_i- a_{ij})^{6}\Big)^{-1} \bigg)^{-1}$\\
Hartmann 3     & $f_{16}(x_1, x_2, x_3) = - \sum\limits_{i = 1}^{4} \alpha_i exp\bigg(-\sum\limits_{j=1}^{3}A_{ij}(x_j - P_{ij})^2\bigg)$\\
\\\vspace{-5mm}\\
\hline
\end{tabular}
\end{table*}

\begin{table*}[htbp]
\centering
\caption{Composite benchmark functions -- adopted from CEC2017 \protect\cite{REF:22}. The definitions of the base functions ($F_N$) are defined in CEC2017.}
\label{tab:composite_benchmarks}
\begin{tabular}{lc}
\hline
\multicolumn{1}{c}{\textbf{Function definition}} & \multicolumn{1}{c}{$f(\vec{x}_{opt})$} \\\hline
\\\vspace{-5mm}\\
$
f_{c,1}(\vec{x}) = \left\{
        \begin{array}{l}
            F_1 = \text{Rotated and Shifted Rosenbrock’s Function}, \\
            F_2 = \text{Rotated and Shifted High Conditioned Elliptic Function},\\
            F_3 = \text{Rotated and Shifted Rastrigin’s Function},\\
            \sigma = [10,20,30], \lambda = [1,\text{$1$E+$06$},1], bias = [0,100,200].
        \end{array}
    \right.
$
\newline

& $2100$
\\~\\
$
f_{c,2}(\vec{x}) = \left\{
        \begin{array}{l}
            F_1 = \text{Rotated and Shifted Rastrigin’s Function}, \\
            F_2 = \text{Rotated and Shifted Griewank’s Function},\\
            F_3 = \text{= Rotated and Shifted Modified Schwefel’s Function},\\
            \sigma = [10,20,30], \lambda = [1,10,1], bias = [0,100,200].
        \end{array}
    \right.
$
& $2200$
\\~\\
$
f_{c,3}(\vec{x}) = \left\{
        \begin{array}{l}
            F_1 = \text{Rotated and Shifted Rosenbrock's Function}, \\
            F_2 = \text{Rotated and Shifted Ackley's Function},\\
            F_3 = \text{Rotated and Shifted Modified Schwefel's Function},\\
            F_4 = \text{Rotated and Shifted Rastrigin's Function},\\
            \sigma = [10,20,30,40], \lambda = [1,10,1,1], bias = [0,100,200,300].
        \end{array}
    \right.
$
& $2300$
\\~\\
$
f_{c,4}(\vec{x}) = \left\{
        \begin{array}{l}
            F_1 = \text{Rotated and Shifted Ackley’s Function}, \\
            F_2 = \text{Rotated and Shifted High Conditioned Elliptic Function},\\
            F_3 = \text{Rotated and Shifted Griewank’s Function},\\
            F_4 = \text{Rotated and Shifted Rastrigin's Function},\\
            \sigma = [10,20,30,40], \lambda = [1,\text{$1$E+$06$},10,1], bias = [0,100,200,400].
        \end{array}
    \right.
$
& $2400$
\\~\\
$
f_{c,5}(\vec{x}) = \left\{
        \begin{array}{l}
            F_1 = \text{Rotated and Shifted Rastirigin’s Function}, \\
            F_2 = \text{Rotated and Shifted HappyCat Function},\\
            F_3 = \text{Rotated and Shifted Ackley’s Function},\\
            F_4 = \text{Rotated and Shifted Discus Function},\\
            F_5 = \text{Rotated and Shifted Rosenbrock's Function},\\
            \sigma = [10,20,30,40,50], \lambda = [10,1,10,\text{$1$E+$06$},1],\\ bias = [0,100,200,300,400].
        \end{array}
    \right.
$
& $2500$
\\~\\
$
f_{c,6}(\vec{x}) = \left\{
        \begin{array}{l}
            F_1 = \text{Rotated and Shifted Expanded Scaffer’s Function}, \\
            F_2 = \text{Rotated and Shifted Modified Schwefel’s Function},\\
            F_3 = \text{Rotated and Shifted Griewank’s Function},\\
            F_4 = \text{Rotated and Shifted Rosenborck's Function},\\
            F_5 = \text{Rotated and Shifted Rastrigin's Function},\\
            \sigma = [10,20,30,40], \lambda = [\text{$1$E+$26$},10,\text{$1$E+$06$},10, \text{$5$E+$04$}],\\ bias = [0,100,200,300,400].
        \end{array}
    \right.
$
& $2600$
\\~\\
\hline
\end{tabular}
\end{table*}

\begin{table*}[htbp]
\centering
\ContinuedFloat
\caption{Composite benchmark functions -- adopted from CEC2017 \protect\cite{REF:22}. The definitions of the base functions ($F_N$) are defined in CEC2017.}
\label{tab:composite_benchmarks_b}
\begin{tabular}{lc}
\hline
\multicolumn{1}{c}{\textbf{Function definition}} & \multicolumn{1}{c}{$f_{c}(\vec{x}_{opt})$} \\\hline
\\\vspace{-5mm}\\
$
f_{c,7}(\vec{x}) = \left\{
        \begin{array}{l}
            F_1 = \text{Rotated and Shifted HGBat Function}, \\
            F_2 = \text{Rotated and Shifted Rastrigin’s Function},\\
            F_3 = \text{Rotated and Shifted Modified Schwefel's Function},\\
            F_4 = \text{Rotated and Shifted BentCigar Function},\\
            F_5 = \text{Rotated and Shifted High Conditioned Elliptic Function},\\
            F_6 = \text{Rotated and Shifted Expanded Scaffer's Function},\\
            \sigma = [10,20,30,40,50,60], \lambda = [10,10,2.5,\text{$1$E+$26$},\text{$1$E+$06$},\text{$5$E+$04$}],\\ bias = [0,100,200,300,400,500].
        \end{array}
    \right.
$
& $2700$
\\~\\
$
f_{c,8}(\vec{x}) = \left\{
        \begin{array}{l}
            F_1 = \text{Rotated and Shifted Ackley’s Function}, \\
            F_2 = \text{Rotated and Shifted Griewank’s Function},\\
            F_3 = \text{Rotated and Shifted Discus Function},\\
            F_4 = \text{Rotated and Shifted Rosenbrock’s Function},\\
            F_5 = \text{Rotated and Shifted HappyCat Function},\\
            F_6 = \text{Rotated and Shifted Expanded Scaffer's Function},\\
            \sigma = [10,20,30,40,50,60], \lambda = [10,10,\text{$1$E+$06$},1,1,\text{$5$E+$04$}],\\ bias = [0,100,200,300,400,500].
        \end{array}
    \right.
$
& $2800$
\\~\\
$
f_{c,9}(\vec{x}) = \left\{
        \begin{array}{l}
            F_1 = \text{Hybrid Function 5}, \\
            F_2 = \text{Hybrid Function 8},\\
            F_3 = \text{Hybrid Function 9},\\
            \sigma = [10,30,50], \lambda = [1,1,1], bias = [0,100,200].
        \end{array}
    \right.
$
& $2900$
\\~\\
$
f_{c,10}(\vec{x}) = \left\{
        \begin{array}{l}
            F_1 = \text{Hybrid Function 5}, \\
            F_2 = \text{Hybrid Function 6},\\
            F_3 = \text{Hybrid Function 7},\\
            \sigma = [10,30,50], \lambda = [1,1,1], bias = [0,100,200].
        \end{array}
    \right.
$
& $3000$
\\~\\\hline
\end{tabular}
\end{table*}

\begin{table*}[htbp]
\centering
\caption{Comparison with WOA \protect\cite{REF:6}.}
\label{tab:WOA}
\begin{tabular}{cllc}
\hline
\multicolumn{1}{c}{\textbf{Function}} & \multicolumn{1}{c}{\textbf{Proposed Approach}} & \multicolumn{1}{c}{\textbf{WOA \protect\cite{REF:6}}} & \multicolumn{1}{c}{\textbf{$f(\vec{x}_{opt})$}} \\ \hline
\\\vspace{-5mm}\\
$f_1$ & 0.775955 $\pm$ 0.222759 &  \textbf{1.41E-30 $\pm$ 4.91E-30} & 0\\
$f_4$ & 5.154556 $\pm$ 2.273995 &  \textbf{0.072581 $\pm$ 0.39747} & 0\\
$f_5$ & 0.759658 $\pm$ 0.251091 &  \textbf{1.06E-21 $\pm$ 2.39E-21} & 0\\
$f_6$ & \textbf{26.778816 $\pm$ 3.804910} &  27.86558 $\pm$ 0.763626 & 0\\
$f_{9}$ & 0.809812 $\pm$ 0.088605 & \textbf{0.000289 $\pm$ 0.001586} & 0\\
$f_{10}$ & \textbf{2.230591 $\pm$ 0.810609} & 7.4043 $\pm$ 9.897572 & 0\\
$f_{11}^\dagger$ & \textbf{-12554.89 $\pm$ 33.2842} & -5080.76 $\pm$ 695.7968 & -12569.48\\
$f_{13}$ & \textbf{-1.031611 $\pm$ 3.73E-05} & -1.03163 $\pm$ 4.2E-07 & -1.0316\\
$f_{15}$ & \textbf{0.998004 $\pm$ 5.98E-10} & 2.111973 $\pm$ 2.498594 & 0.998004 $\approx$ 1\\
\\\vspace{-5mm}\\
\hline
\multicolumn{3}{c}{$^\dagger$\scriptsize{This version of Schwefel min $f(\vec{x}_{opt}) =$ -418.9829$\times$30} (see \cite{REF:6})}
\end{tabular}
\end{table*}
% \footnote{}
\begin{table*}[htbp]
\centering
\caption{Comparison with PFA \protect\cite{REF:7}.}
\label{tab:PFA}
\begin{tabular}{cllc}
\hline
\multicolumn{1}{c}{\textbf{Function}} & \multicolumn{1}{c}{\textbf{Proposed Approach}} & \multicolumn{1}{c}{\textbf{PFA \protect\cite{REF:7}}} & \multicolumn{1}{c}{\textbf{$f(\vec{x}_{opt})$}} \\ \hline
\\\vspace{-5mm}\\
$f_2$ & 0.117980 $\pm$ 0.048295 &  \textbf{5.5674E-25 $\pm$ 7.9092E-25} & 0\\
$f_3$ & 0.031421 $\pm$ 0.016130 &  \textbf{9.9813E-46 $\pm$ 3.3585E-45} & 0\\
$f_5$ & 0.437154 $\pm$ 0.249896 &  \textbf{3.4831E-14 $\pm$ 6.2094E-14} & 0\\
$f_7$ & -49.996395 $\pm$ 0.004986 &  \textbf{-50.0000 $\pm$ 1.98E-11} & -50$^\dagger$\\
$f_{8}$ & \textbf{0.081319 $\pm$ 0.027990} & 11.5480 $\pm$ 12.9802 & 0\\
$f_{9}$ & 0.425310 $\pm$ 0.408112 & \textbf{0.0006 $\pm$ 0.0012} & 0\\
$f_{11}$ & \textbf{0.610558 $\pm$ 0.145078} & 3.1549E+3 $\pm$ 5.6274E+2 & 0\\
$f_{14}$ & 3.000043 $\pm$ 7.51E-5 & 3.0000 $\pm$ 2.4952E-16 & 3\\
$f_{16}$ & -3.855772 $\pm$ 0.009925 & \textbf{-3.8628 $\pm$ 1.5026E-15} & -3.8628\\
\\\vspace{-5mm}\\
\hline
\multicolumn{3}{l}{$^\dagger$\scriptsize{This value has been evaluated for $n=6$} (as in \cite{REF:7})}
\end{tabular}
\end{table*}

\begin{table*}[htbp]
\centering
\caption{Scalability with PFA \protect\cite{REF:7}.}
\label{tab:PFAScale}
\begin{tabular}{clllc}
\hline
\textbf{Function} & \textbf{$n$} & \multicolumn{1}{c}{\textbf{Proposed Approach}} & \multicolumn{1}{c}{\textbf{PFA \protect\cite{REF:7}}} & \multicolumn{1}{c}{\textbf{$f(\vec{x}_{opt})$}} \\\hline
\\\vspace{-5mm}\\
\multirow{3}{*}{$f_8$} & 10 & 0.0143 $\pm$ 0.0206 & \textbf{2.17E-35 $\pm$ 1.09E-34} & 0 \\
 & 50 & \textbf{53.6161 $\pm$ 11.3012} & 343.1244 $\pm$ 92.1694 & 0 \\
 & 100 & \textbf{680.8374 $\pm$ 61.9845} & 1911.1759 $\pm$ 195.2386 & 0 \\
\hline \\\vspace{-5mm}\\
\multirow{3}{*}{$f_{11}$} & 10 & \textbf{0.3395 $\pm$ 0.91820} & 449.5090 $\pm$ 184.2821 & 0 \\
 & 50 & \textbf{129.6466 $\pm$ 82.8948} & 6423.8391 $\pm$ 663.2207 & 0 \\
 & 100 & \textbf{7208.6969 $\pm$ 489.7975} & 13610.2227 $\pm$ 1361.9760 & 0 \\
 \\\vspace{-5mm}\\
 \hline
\end{tabular}
\end{table*}

\begin{table*}[htbp]
\centering
\caption{Comparison with MPSFHS \protect\cite{REF:1a}.}
\label{tab:MPSFHS}
\begin{tabular}{cllc}
\hline
\multicolumn{1}{c}{\textbf{Function}} & \multicolumn{1}{c}{\textbf{Proposed Approach}} & \multicolumn{1}{c}{\textbf{MPSFHS \protect\cite{REF:1a}}} & \multicolumn{1}{c}{\textbf{$f(\vec{x}_{opt})$}} \\ \hline
\\\vspace{-5mm}\\
$f_{8}$ & \textbf{173.1648 $\pm$ 17.1974} &  1181.5644 $\pm$ 71.3395 & 0 \\
$f_{11}$ & \textbf{971.8366 $\pm$ 150.2988} &  1554.8070 $\pm$ 178.5386 & 0 \\
\\\vspace{-5mm}\\
\hline
\end{tabular}
\end{table*}

% \begin{table*}[htbp]
{\captionsetup[table]{aboveskip=2cm}\begin{sidewaystable*}

\centering
\caption{CEC2017 composite functions tested with $n=2$ by using only $10$ iterations (only for the PSS, PFA, WOA, and PSO and equivalent evaluations for the MPSFHS) for $f_{c,1} \to f_{c,8}$ -- reported the error values here $\big(f_{c,i}(\vec{x_{i}}) - f_{c,i}(\vec{x}_{opt})\big)$.}
\label{tab:PSS_COMPO}
\begin{tabular}{lllcccccccc}
\hline
\textbf{Algorithm} & & $f_{c,1}$ & $f_{c,2}$ & $f_{c,3}$ & $f_{c,4}$ & $f_{c,5}$ & $f_{c,6}$ & $f_{c,7}$ & $f_{c,8}$ \\ \hline
\\\vspace{-5mm}\\
\multirow{5}{*}{PSS} & min      & 2.4285E-02 & 6.5834E-02 & 6.0717E-02 & 7.2868E+00 & 2.8464E-01 & 6.2283E-02 & 6.4629E-01 & 5.5129E+00\\
& max      & 1.0018E+02 & 1.1226E+02 & 3.0392E+02 & 2.0635E+02 & 2.1118E+02 & 2.0057E+02 & 1.8914E+02 & 2.1992E+02\\
& median   & 4.1376E+00 & 1.2610E+01 & 7.0064E+01 & 1.0059E+02 & 1.0614E+02 & 1.7070E+00 & 2.6410E+00 & 1.0117E+02\\
& mean     & 1.9017E+01 & 2.5048E+01 & 9.4231E+01 & 1.0793E+02 & 1.1002E+02 & 2.7532E+01 & 2.7133E+01 & 1.0510E+02\\
& std      & 3.1468E+01 & 3.3741E+01 & 1.1463E+02 & 5.5402E+01 & 6.1274E+01 & 5.2681E+01 & 5.0566E+01 & 7.5755E+01\\ \\\vspace{-5mm}\\\hline\\ 

\multirow{5}{*}{PFA} & min      & 5.1962E-02 & 2.0261E-02 & 7.8146E-01 & 2.1461E+01 & 6.7842E+00 & 8.8438E-02 & 1.3907E+00 & 1.2169E+01 \\
& max      & 1.0021E+02 & 1.0230E+02 & 3.0303E+02 & 2.0244E+02 & 3.0168E+02 & 2.0858E+02 & 4.0021E+02 & 3.5366E+02 \\
& median   & 4.0589E-01 & 8.0329E-01 & 9.6471E+00 & 1.0142E+02 & 1.0713E+02 & 1.9240E+00 & 3.6875E+00 & 1.1203E+02 \\
& mean     & 4.2739E+00 & 9.5696E+00 & 9.5648E+01 & 1.0874E+02  & 1.1802E+02 & 9.1956E+00 & 1.7240E+01 & 1.3732E+02 \\
& std      & 1.8221E+01 & 2.5488E+01 & 1.3689E+02 & 3.7036E+01 & 8.2871E+01 & 3.7697E+01 & 7.2359E+01 & 7.4851E+01 \\ \\\vspace{-5mm}\\\hline\\

\multirow{5}{*}{WOA} & min      & 3.9120E-02 & 4.6032E-02 & 3.8543E-01 & 1.4056E+01 & 3.4195E+00 & 3.7880E-01 & 2.4574E+00 & 3.2645E+01 \\
& max      & 1.1761E+02 & 2.3130E+01 & 3.0207E+02 & 2.1358E+02 & 3.0999E+02 & 2.0312E+02 & 1.8917E+02 & 3.0409E+02 \\
& median   & 6.1233E-01 & 1.0885E+01 & 5.3677E+00 & 1.0070E+02 & 1.1937E+02 & 5.3988E+00 & 2.4332E+01 & 1.6472E+02 \\
& mean     & 4.9854E+00 & 8.9399E+00 & 1.0609E+02 & 1.0771E+02 & 1.3170E+02 & 2.2009E+01 & 4.8168E+01 & 1.6502E+02 \\
& std      & 2.1319E+01 & 5.9659E+00 & 1.2264E+02 & 3.3913E+01 & 9.4636E+01 & 4.6116E+01 & 5.4586E+01 & 6.5213E+01 \\ \\\vspace{-5mm}\\\hline\\ 

\multirow{5}{*}{MPSFHS} & min      & 1.3582E+00 & 3.0120E-01 & 8.8972E+00 & 9.2268E+01 & 3.1234E+01 & 1.6434E-01 & 3.5703E+00 & 5.2473E+01 \\
& max      & 6.0320E+01 & 1.1324E+02 & 3.0430E+02 & 2.0983E+02 & 3.1993E+02 & 1.4941E+02 & 7.8858E+01 & 3.0867E+02 \\
& median   & 1.9290E+01 & 2.4465E+01 & 7.5047E+01 & 1.3201E+02 & 2.0406E+02 & 4.1046E+01 & 2.8917E+01 & 1.4537E+02 \\
& mean     & 2.1165E+01 & 3.0313E+01 & 1.0591E+02 & 1.4425E+02 & 1.8885E+02 & 4.8282E+01 & 3.3665E+01 & 1.5904E+02 \\
& std      & 1.7733E+01 & 2.3374E+01 & 8.9850E+01 & 3.8189E+01 & 7.9865E+01 & 3.5734E+01 & 1.9670E+01 & 5.5776E+01 \\ \\\vspace{-5mm}\\\hline\\ 

\multirow{5}{*}{PSO} & min      & 3.0762E-01 & 6.8406E-01 & 2.4319E+00 & 7.1158E+01 & 1.1524E+01 & 2.9869E+00 & 5.9283E+00 & 4.4153E+01 \\
& max      & 1.0109E+02 & 1.7330E+01 & 3.0334E+02 & 1.6073E+02 & 2.8149E+02 & 2.2489E+02 & 4.0218E+02 & 2.1841E+02 \\
& median   & 5.8090E+00 & 7.5109E+00 & 2.6452E+01 & 1.1033E+02 & 1.4209E+02 & 2.0422E+01 & 2.1516E+01 & 1.2396E+02 \\
& mean     & 9.4914E+00 & 8.2878E+00 & 1.1008E+02 & 1.1224E+02 & 1.4144E+02 & 2.7641E+01 & 3.4890E+01 & 1.3874E+02 \\
& std      & 1.7832E+01 & 4.8397E+00 & 1.3079E+02 & 1.5796E+01 & 7.7616E+01 & 3.8640E+01 & 7.0202E+01 & 4.8416E+01 \\ \\\vspace{-5mm}
% \\\hline\\ 

% \multirow{5}{*}{CS} & min      &  &  &  &  &  &  &  &  \\
% & max      &  &  &  &  &  &  &  &  \\
% & median   &  &  &  &  &  &  &  &  \\
% & mean     &  &  &  &  &  &  &  &  \\
% & std      &  &  &  &  &  &  &  &  \\
\\\vspace{-5mm}\\
\hline
\end{tabular}
\end{sidewaystable*}
}%

\begin{table*}[htbp]

\centering
\caption{CEC2017 composite functions tested with $n=100$ by using only $100$ iterations for both the PSS and PFA algorithms on $f_{c,1} \to f_{c,10}$ -- reported the error values here $\big(f_{c,i}(\vec{x_{i}}) - f_{c,i}(\vec{x}_{opt})\big)$.}
\label{tab:PSS_COMPO2}
\begin{tabular}{lllccccc}
\hline
\textbf{Algorithm} & & $f_{c,1}$ & $f_{c,2}$ & $f_{c,3}$ & $f_{c,4}$ & $f_{c,5}$ \\ \hline
\\\vspace{-5mm}\\
\multirow{13}{*}{PSS} & min      & 1.2002E+03 & 2.5724E+04 & 1.5424E+03 & 2.1423E+03 & 5.1997E+03 \\
& max      & 1.7306E+03 & 3.1503E+04 & 1.8003E+03 & 2.5840E+03 & 9.2997E+03 \\
& median   & 1.5094E+03 & 2.8014E+04 & 1.6500E+03 & 2.3455E+03 & 7.0022E+03 \\
& mean     & 1.4929E+03 & 2.7920E+04 & 1.6588E+03 & 2.3551E+03 & 7.0311E+03 \\
& std      & 1.0472E+02 & 1.3480E+03 & 7.0691E+01 & 1.1244E+02 & 1.0787E+03 \\ \\\cline{2-8}
 & & $f_{c,6}$ & $f_{c,7}$ & $f_{c,8}$ & $f_{c,9}$ & $f_{c,10}$\\\cline{2-8}\\
& min      & 1.6692E+04 & 1.5568E+03 & 4.9472E+03 & 5.1738E+03 & 1.7476E+07 \\
& max      & 1.9222E+04 & 2.2142E+03 & 1.2336E+04 & 7.4074E+03 & 8.2332E+07 \\
& median   & 1.7630E+04 & 1.6990E+03 & 8.4480E+03 & 6.3363E+03 & 4.2002E+07 \\
& mean     & 1.7686E+04 & 1.7340E+03 & 8.7575E+03 & 6.2897E+03 & 4.2233E+07 \\
& std      & 8.4035E+02 & 1.4641E+02 & 1.9430E+03 & 5.6462E+02 & 1.6223E+07 \\ \\\hline
\textbf{Algorithm} & & $f_{c,1}$ & $f_{c,2}$ & $f_{c,3}$ & $f_{c,4}$ & $f_{c,5}$ \\ \hline \\
\multirow{13}{*}{PFA} & min      & 1.3493E+03 & 2.4074E+04 & 2.3835E+04 & 2.1928E+03 & 5.6223E+03 & \\
& max      & 1.9686E+03 & 3.2000E+04 & 3.3604E+04 & 2.8422E+03 & 1.6262E+04 \\
& median   & 1.5668E+03 & 2.6900E+04 & 2.8106E+04 & 2.5436E+03 & 6.9669E+03 \\
& mean     & 1.6031E+03 & 2.7021E+04 & 2.7941E+04 & 2.5141E+03 & 7.7530E+03 \\
& std      & 1.5386E+02 & 2.2703E+03 & 2.7965E+03 & 1.7907E+02 & 2.2391E+03 \\  \\\cline{2-8}
 & & $f_{c,6}$ & $f_{c,7}$ & $f_{c,8}$ & $f_{c,9}$ & $f_{c,10}$\\\cline{2-8}\\
& min      & 1.6504E+04 & 1.1039E+03 & 7.0900E+03 & 5.7046E+03 & 4.2582E+07 \\
& max      & 2.2592E+04 & 1.8703E+03 & 1.6716E+04 & 1.0894E+04 & 2.0416E+08 \\
& median   & 1.9795E+04 & 1.4006E+03 & 1.1057E+04 & 7.7024E+03 & 9.5209E+07 \\
& mean     & 1.9688E+04 & 1.4345E+03 & 1.1339E+04 & 7.8995E+03 & 1.0203E+08 \\
& std      & 1.5073E+03 & 1.8877E+02 & 2.3745E+03 & 1.2773E+03 & 4.5804E+07 \\ \\\hline
\end{tabular}
\end{table*}
%  \\\vspace{-5mm}\\\hline\\ 
% \multirow{5}{*}{CS} & min      &  &  &  &  &  &  &  &  \\
% & max      &  &  &  &  &  &  &  &  \\
% & median   &  &  &  &  &  &  &  &  \\
% & mean     &  &  &  &  &  &  &  &  \\
% & std      &  &  &  &  &  &  &  &  \\
% \end{table*}

{\captionsetup[table]{aboveskip=2cm}\begin{sidewaystable*}
\centering
\caption{Composite functions comparison with PFA algorithm \protect\cite{REF:7} -- reported here the error values $\big(f_{c,i}(\vec{x_{i}}) - f_{c,i}(\vec{x}_{opt})\big)$.}
% \Rotatebox{90}{%
\resizebox{1.00\columnwidth}{!}{%
\centering
\label{tab:compbench}
\begin{tabular}{lclcccccccccc}
\hline
\textbf{Method}                & $n$                   & \textbf{Error}     & $f_{c,1}$     & $f_{c,2}$     & $f_{c,3}$     & $f_{c,4}$     & $f_{c,5}$     & $f_{c,6}$     & $f_{c,7}$     & $f_{c,8}$    & $f_{c,9}$     & $f_{c,10}$ \\ \hline
\\\vspace{-5mm}\\
\multirow{15}{*}{PFA \protect\cite{REF:7}} & \multirow{5}{*}{10} & min         &   1.0000E+02  &   1.1563E+01  &   3.0014E+02  &   1.0006E+02  &   3.9774E+02  &   3.0000E+02  &   3.8848E+02  &   3.0000E+02  &   2.4124E+02  & 2.8612E+03     \\
                      &                     & max         &   2.2019E+02  &   1.0311E+02  &   3.1846E+02  &   3.4380E+02  &   4.4602E+02  &   3.0000E+02  &   3.9570E+02  &   4.4498E+02  &   3.2367E+02  & 9.2812E+04     \\
                      &                     & median      &   1.5467E+02  &   1.0137E+02  &   3.0762E+02  &   1.6547E+02  &   4.0078E+02  &   3.0000E+02  &   3.9090E+02  &   3.0277E+02  &   2.8414E+02  & 1.7656E+04     \\
                      &                     & mean        &   1.5645E+02  &   9.3593E+01  &   3.0706E+02  &   2.1330E+02  &   4.1863E+02  &   3.0000E+02  &   3.9175E+02  &   3.4087E+02  &   2.8418E+02  &  6.2390E+04    \\
                      &                     & std         &   5.5756E+01  &   2.4340E+01  &   4.3861E+00  &   9.8189E+01  &   2.1735E+01  &   2.1902E-05  &   1.9096E+00  &   6.3004E+01  &   1.5359E+01  &  1.5105E+04    \\ \cline{2-13} 
                      & \multirow{5}{*}{30} & min         &   2.0778E+02  &   1.0000E+02  &   3.7323E+02  &   5.2936E+02  &   3.8300E+02  &   2.0782E+02  &   4.8804E+02  &   4.0528E+02  &   5.6256E+02  &  2.4685E+03    \\
                      &                     & max         &   4.2379E+02  &   7.5775E+03  &   5.2140E+02  &   5.8619E+02  &   4.3544E+02  &   2.6365E+03  &   5.1048E+02  &   4.7741E+02  &   1.0563E+02  &  1.5743E+05    \\
                      &                     & median      &   3.1604E+02  &   1.9984E+03  &   4.9308E+02  &   5.7169E+02  &   3.8810E+02  &   9.5600E+02  &   5.0277E+02  &   4.3275E+02  &   6.8022E+02  &  5.7236E+04    \\
                      &                     & mean        &   3.2029E+02  &   3.3880E+03  &   4.8791E+02  &   5.7049E+02  &   3.8989E+02  &   9.5505E+02  &   5.0281E+02  &   4.3839E+02  &   7.1405E+02  &  6.4199E+04    \\
                      &                     & std         &   5.5865E+01  &   3.3643E+03  &   2.6840E+01  &   1.2125E+01  &   7.3578E+00  &   6.2841E+02  &   2.5284E+00  &   2.5585E+01  &   1.1731E+02  &  3.0136E+04    \\ \cline{2-13} 
                      & \multirow{5}{*}{50} & min         &   2.9169E+02  &   1.0000E+02  &   5.1701E+02  &   5.5598E+02  &   5.4014E+02  &   3.5236E+02  &   5.0000E+02  &   5.1385E+02  &   4.6940E+02  &  1.1016E+06    \\
                      &                     & max         &   6.5178E+02  &   1.4177E+04  &   7.7342E+02  &   8.3818E+02  &   6.7421E+02  &   4.1024E+03  &   6.8090E+02  &   6.5809E+02  &   1.4784E+03  &  1.8261E+06    \\
                      &                     & median      &   5.6687E+02  &   1.3255E+04  &   7.2031E+02  &   8.1295E+02  &   6.0900E+02  &   2.0229E+03  &   6.2269E+02  &   5.5260E+02  &   8.7184E+02  &  1.4303E+06    \\
                      &                     & mean        &   5.1824E+02  &   1.3002E+04  &   7.0909E+02  &   7.9178E+02  &   6.1075E+02  &   2.2147E+03  &   6.2270E+02  &   5.6326E+02  &   9.0843E+02  &  1.4301E+06    \\
                      &                     & std         &   1.0989E+01  &   1.9255E+03  &   4.9710E+01  &   6.0605E+01  &   2.1206E+01  &   6.3598E+02  &   3.8000E+01  &   3.0291E+01  &   2.5833E+02  &  1.5398E+05   \\
                      \\\vspace{-5mm}\\ \cline{1-13} \\\vspace{-5mm}\\

\multirow{15}{*}{PSS} & \multirow{5}{*}{10} & min         &   1.0000E+02   &   2.5811E+01   &                                  3.0543E+02   &   1.0004E+02   &   1.0549E+02   &   5.0521E+00   &                          3.8891E+02   &   3.9683E+01   &   2.3791E+02   &    2.7732E+03   \\
                      &                     & max         &   2.6412E+02   &   1.1087E+02   &   3.4689E+02   &   3.8552E+02   &   4.4623E+02   &   4.5287E+02   &   4.4213E+02   &   6.1188E+02   &   3.1690E+02   &    8.4283E+04   \\
                      &                     & median      &   2.2392E+02   &   1.0795E+02   &   3.2332E+02   &   3.5939E+02   &   3.9910E+02   &   3.0109E+02   &   4.0022E+02   &   3.0352E+02   &   2.6677E+02   &    1.3233E+04   \\
                      &                     & mean        &   2.0725E+02   &   1.0058E+02   &   3.2451E+02   &   3.1612E+02   &   3.9322E+02   &   2.9345E+02   &   4.0051E+02   &   3.4762E+02   &   2.7031E+02   &    1.9025E+04   \\
                      &                     & std         &   4.8079E+01   &   1.9823E+01   &   1.0129E+01   &   1.0004E+02   &   5.2261E+01   &   5.5000E+01   &   8.1409E+00   &   1.1675E+02   &   1.7654E+01   &    1.6283E+04   \\ \cline{2-13} 
                      & \multirow{5}{*}{30} & min         &   2.5982E+02   &   9.5178E+01   &   4.1179E+02   &   5.0458E+02   &   3.8376E+02   &   2.1092E+02   &   4.9558E+02   &   4.0002E+02   &   4.7401E+02   &    1.2650E+04   \\
                      &                     & max         &   3.5442E+02   &   3.0117E+03   &   5.0432E+02   &   6.5072E+02   &   3.9624E+02   &   2.2463E+03   &   5.5834E+02   &   4.6109E+02   &   1.0332E+03   &    1.8496E+06   \\
                      &                     & median      &   2.9840E+02   &   1.0773E+02   &   4.6186E+02   &   5.9202E+02   &   3.8751E+02   &   3.0266E+02   &   5.2653E+02   &   4.1788E+02   &   6.9333E+02   &    6.0821E+04   \\
                      &                     & mean        &   2.9769E+02   &   5.9879E+02   &   4.6472E+02   &   5.8596E+02   &   3.8745E+02   &   7.1025E+02   &   5.2544E+02   &   4.2076E+02   &   6.8941E+02   &    1.0582E+05   \\
                      &                     & std         &   1.6933E+01   &   9.9855E+02   &   2.2548E+01   &   3.1120E+01   &   1.8791E+00   &   7.4929E+02   &   1.2494E+01   &   1.2597E+01   &   1.1408E+02   &    2.5646E+05   \\ \cline{2-13} 
                      & \multirow{5}{*}{50} & min         &   3.2316E+02   &   1.2346E+02   &   5.8261E+02   &   7.1224E+02   &   5.1911E+02   &   3.1799E+02   &   5.8963E+02   &   4.6972E+02   &   6.8860E+02   &    1.8473E+06   \\
                      &                     & max         &   4.2892E+02   &   6.5096E+03   &   7.8839E+02   &   9.0554E+02   &   6.2009E+02   &   3.9818E+03   &   8.0712E+02   &   6.1470E+02   &   1.9981E+03   &    3.8191E+06   \\
                      &                     & median      &   3.8511E+02   &   5.4346E+03   &   6.7549E+02   &   8.4179E+02   &   5.7757E+02   &   3.4853E+03   &   6.9824E+02   &   5.2609E+02   &   1.1805E+03   &    2.5966E+06 
  \\
                      &                     & mean        &   3.8574E+02   &   5.3759E+03   &   6.7342E+02   &   8.3583E+02   &   5.7406E+02   &   3.2646E+03   &   7.0768E+02   &   5.2915E+02   &   1.2065E+03   &   2.5769E+06    \\
                      &                     & std         &   2.3097E+01   &   8.8643E+02   &   3.6980E+01   &   4.3656E+01   &   2.6876E+01   &   7.4941E+02   &   5.3367E+01   &   2.3379E+01   &   2.4202E+02   &   3.5785E+05   
                      \\\vspace{0mm}\\ \cline{1-13}

\end{tabular}
% %
% }
}%
\end{sidewaystable*}} % <- closing brace

% Graphical results for composite functions
\begin{figure*}[htbp]
\centerline{\includegraphics[width = 1.0\linewidth]{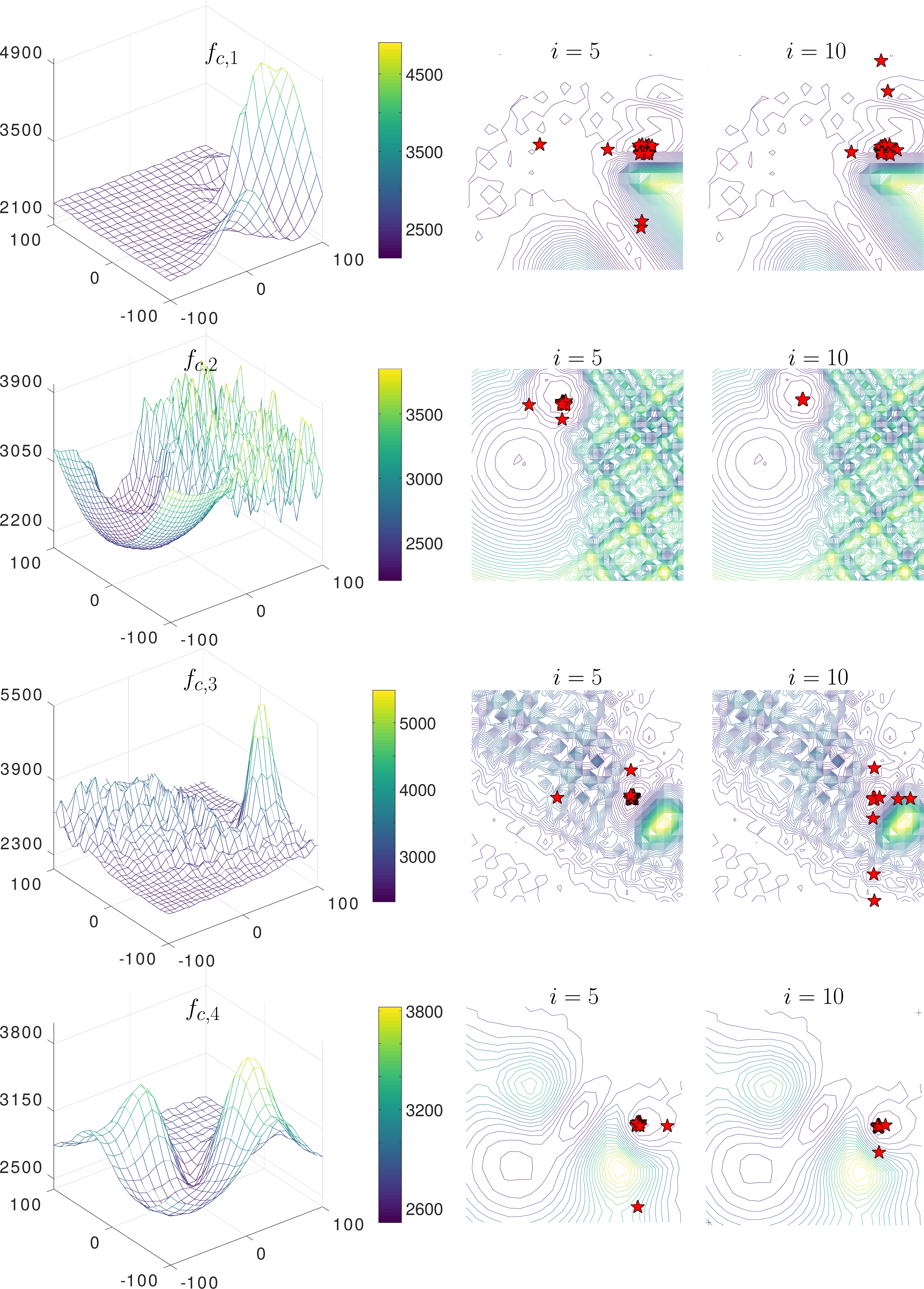}}
\caption{The convergence of composite functions using $\alpha = 0.95$, $\beta = 30$, and $\gamma = 10$ -- illustrated the best runs for $f_{c,1} \to f_{c,4}$ where $n=2$.}
\label{FIG:9}
\end{figure*}

\begin{figure*}[htbp]
\ContinuedFloat
\centerline{\includegraphics[width = 1.0\linewidth]{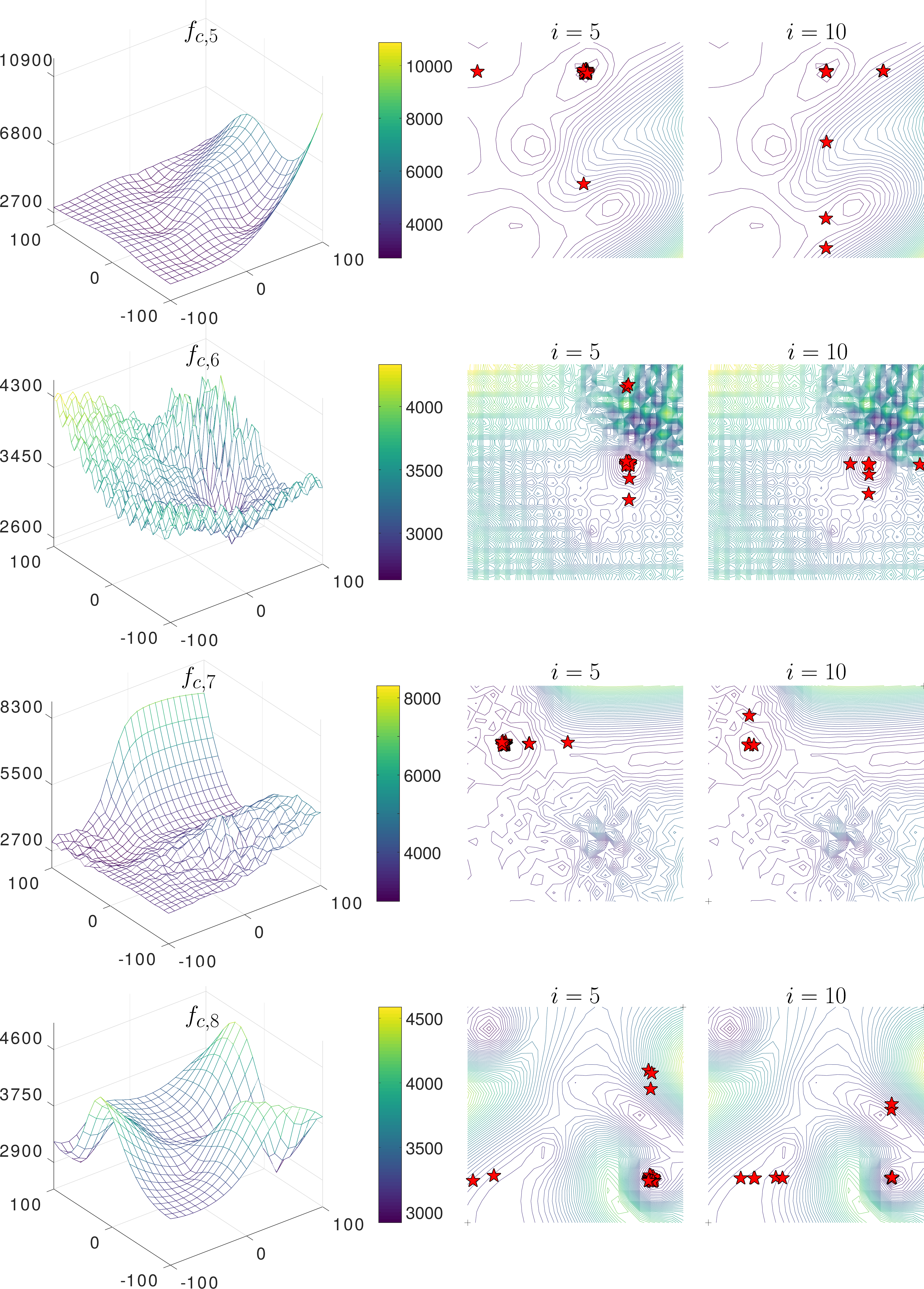}}
\caption{The convergence of composite functions using $\alpha = 0.95$, $\beta = 30$, and $\gamma = 10$ -- illustrated the best runs for $f_{c,5} \to f_{c,8}$ where $n=2$.}
\label{FIG:9b}
\end{figure*}

\subsection{Engineering Benchmarks}
\label{subsec:engbenchmarking}
\subsubsection{Traditional Engineering Benchmarks}
\begin{figure*}[htbp]
\centerline{\includegraphics[width = 1.2\linewidth]{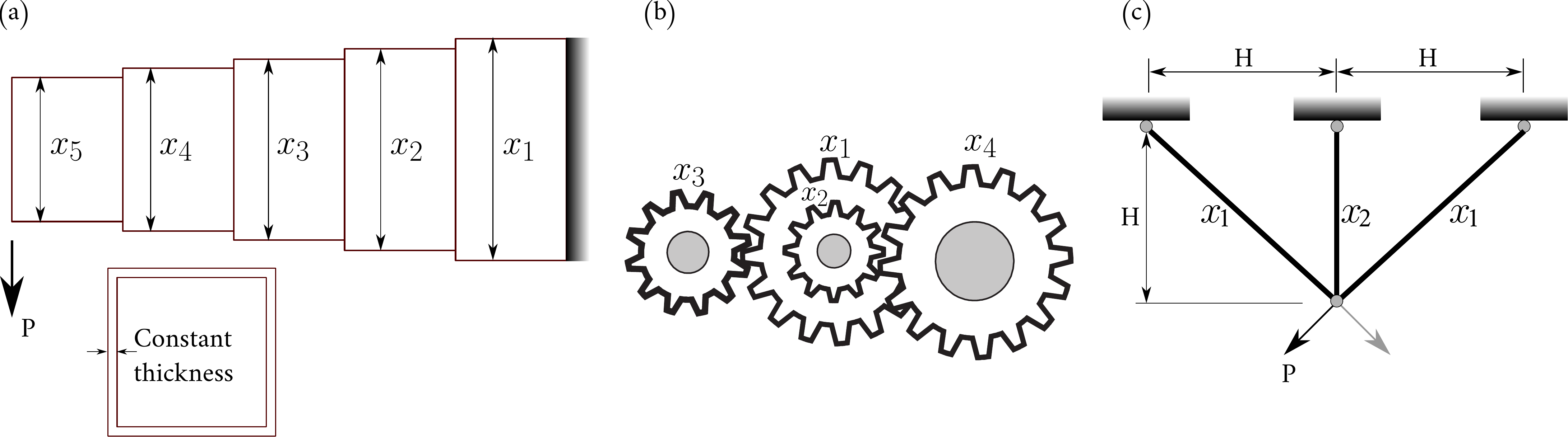}}
\caption{(a) The cantilever beam design problem; (b) the train gears design problem (a modified figure from \protect\cite{REF:8}); (c) the three-bar truss design problem.}
\label{FIG:6}
\label{FIG:7}
\label{FIG:8}
\end{figure*}

In this section, the PSS was tested and compared with some of the traditional constrained engineering design problems, including: \rom{1}) the cantilever beam design (Fig. \ref{FIG:6} (a)), \rom{2}) the train gears design (Fig. \ref{FIG:7} (b)), and \rom{3}) the three-bar truss design (Fig. \ref{FIG:8} (c)). The reader is referred to Yapici and Cetinkaya (2019) \cite{REF:7} and Mirjalili (2015) \cite{REF:8} for the formulations of these problems.

% \begin{figure}[htbp]
% \centerline{\includegraphics[width = 0.45\linewidth]{cantilever_beam.pdf}}
% \caption{The cantilever beam design problem}
% \label{FIG:6}
% \end{figure}

% \begin{figure}[htbp]
% \centerline{\includegraphics[width = 0.45\linewidth]{train_gear.pdf}}
% \caption{The train gears design problem--a modified figure from \protect\cite{REF:8}}
% \label{FIG:7}
% \end{figure}

% \begin{figure}[htbp]
% \centerline{\includegraphics[width = 0.45\linewidth]{three-bars-truss.pdf}}
% \caption{The three-bar truss design problem}
% \label{FIG:8}
% \end{figure}

The results of the three problems are shown in Table \ref{tab:engineering}. For the first problem, we obtained an almost-identical solution to the PFA algorithm. However, Sharma and Sah recently obtained an even better optimum solution of $1.32545$ with the novel hybrid algorithm m-MBOA Sharma and Saha (2019) \cite{REF:10}. Another recently published article obtained an optimum solution of $1.330565414$ by \cite{REF:56}.

Regarding the train gears design problem, the PSS algorithm reached the same optimal answer as many other algorithms using integer formulation, such as the Moth-flame Optimization Algorithm (MOA) Mirjalili (2015) \cite{REF:8}. Notwithstanding, Sharma and Saha (2019) \cite{REF:10} claim to obtain the optimum answer of $3.3610E-16$ by considering the problem landscape as continuous. By also considering the problem's distribution to be continuous, we obtained a solution of $2.2695E-21$ with $\vec{x} =$ $\{43.170946,$ $14.205443,$ $ 17.919979,$ $40.869247\}$. However, these solution must not be considered as a new valid optimum, as it disregards the nature of the problem by dealing with the number of teeth per gear as a continuously distributed design variable.

Finally, the third case study minimizes the weight of a three-bar truss structure. As can be seen in Table \ref{tab:engineering}, the PSS obtained a slightly better solution than the one obtained by \cite{REF:8}. Sharma and Saha (2019) \cite{REF:10} again detail a new optimal solution of $1.325454889710144$, which lays outside of the feasibility domain of the problem and is therefore inconsequential for engineers. We believe that this solution is obtained by incorrectly enforcing the constraints, or it could be a simple typo. It is worth to mention that the recent Water Strider Algorithm (WSA) by Kaveh and Eslamlou (2020) \cite{REF:32} obtained a slightly better solution with a weight of $263.89584340$ kN.

\subsubsection{Design of Reinforced Concrete Beams}
\begin{figure*}[htbp]
\centerline{\includegraphics[width = 1.5\linewidth]{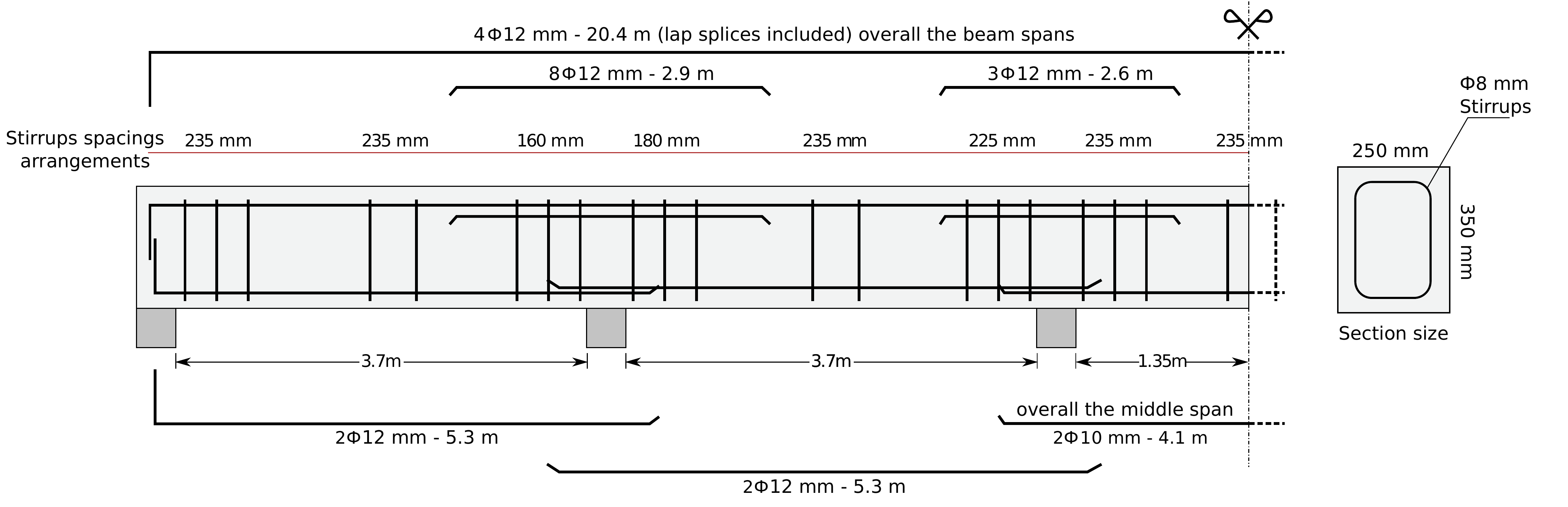}}
\caption{Case study (1) retrieved from \protect\cite{REF:1a} -- best reinforcement details by PSS.}
\label{FIG:2}
\end{figure*}

In the fourth design case study, we recall the problem of designing a reinforced concrete beam (formulated in Shaqfa and Orb\'an (2019) \cite{REF:1a}). In this example, the PSS algorithm solved a complex weight minimization case study that was presented and solved by \cite{REF:1a} using the MPSFHS algorithm. Herein, we used an $\alpha$ of $0.85$ instead of $0.95$ as previously, because this problem required more diversification. Table \ref{tab:engineering} shows the difference between the results obtained with both the MPSFHS and the PSS algorithms. For this discrete problem with $25$ independent design variables, the best answer was obtained by the MPSFHS. However, the difference between the two is small and can be related to some details of the top reinforcement (see Fig. \ref{FIG:2}). This level of perturbation normally requires a small-stepped random walk algorithm to exploit small details. As outlined above, the PSS algorithm does not perform well with regard to finding the local optimum and should be, for this purpose, combined with other algorithms whenever necessary.

\begin{table}[htbp]
\centering
\caption{Comparisons of solutions to engineering problems.}
\label{tab:engineering}
\begin{tabular}{ccc}
\multicolumn{3}{c}{\textbf{\rom{1}) Cantilever Beam Design}}\\ \hline
\textbf{Method} & \textbf{Proposed Approach} & \textbf{PFA \cite{REF:7}}\\ \hline
Weight & 1.33995664399519 & \textbf{1.33995638}\\
$x_1$ & 6.01683010096092  & 6.0154633\\
$x_2$ & 5.30655187659779  & 5.3090222\\
$x_3$ & 4.49420948422588  & 4.4946314\\
$x_4$ & 3.50272928517748  & 3.5017850\\
$x_5$ & 2.15334341962752  & 2.1527578\\
\\\vspace{-5mm}\\
\multicolumn{3}{c}{\textbf{\rom{2}) Train Gears Design}}\\ \hline
\textbf{Method} & \textbf{Proposed Approach} & \textbf{MFO \cite{REF:8}}\\ \hline
Gear ratio  & 2.7009E-12 & 2.7009E-12\\
$x_1$ & 43 & 43\\
$x_2$ & 19 & 19\\
$x_3$ & 16 & 16\\
$x_4$ & 49 & 49\\
\\\vspace{-5mm}\\
\multicolumn{3}{c}{\textbf{\rom{3}) Three-bar Truss Design}}\\ \hline
\textbf{Method} & \textbf{Proposed Approach} & \textbf{MFO \cite{REF:8}}\\ \hline
Weight      & \textbf{263.895843501333} & 263.895979682\\
$x_1$       & 0.788683438026281         & 0.788244770931922\\
$x_2$       & 0.408224806061712         & 0.409466905784741\\
\\\vspace{-5mm}\\
\multicolumn{3}{c}{\textbf{\rom{4}) Reinforced Concrete Beam Design}}\\ \hline
\textbf{Method} & \textbf{Proposed Approach} & \textbf{MPSFHS \cite{REF:1a}}\\ \hline
Weight [kN] & 42.138 & \textbf{42.125}\\
\hline
\end{tabular}
\end{table}

% ********* Conclusions *********
\section{Conclusions}
\label{sec:conclusions}
In this work, we propose a heuristic approach that uses a simple analogy with classical DOE methods. The presented approach samples solution features more densely in a detected prominent region than in the entire search domain. The design of this algorithm was kept as simple as possible, while avoiding structural bias and premature convergence. The algorithm requires three input parameters, namely the population size ($\beta$), the maximum number of iterations ($\gamma$), and the acceptance probability ($\alpha$). We provided a simple probabilistic derivation for investigating good parameter choices. The algorithm was benchmarked against state-of-the-art algorithms using a selected set of standard and engineering optimization problems. The algorithm behaved better in allocating global minima for nonconvex and multimodal problems than the state-of-the-art algorithms, while it does need enhancements on the intensification side to make better use of all the available solution vectors. Moreover, the PSS algorithm proved useful and outperformed state-of-the-art algorithmic proposals in the scalability problems and under extreme cases where only a few iterations are allowed. Even though we used it as a standalone metaheuristic in this paper when benchmarking the algorithm, the algorithm will likely be most useful hybridised in sequential or parallel schemes with other well-known algorithms, e.g., L\'evy flights.

% ********* Future work *********
\section{Future work}
\label{sec:future}
In future work, we plan to implement the \textit{lbest} (see Eberhart and Kennedy (1995) \cite{REF:14}) analogy with the algorithm, where it can redefine the prominent domain region, $\Omega^{'}$, by using, for instance, the best $(1-\alpha) \times \beta$ candidates per iteration. This topology could be used to strengthen the intensification of the algorithm. Hybridization via other well-known algorithms, including but not limited to L\'evy flights \cite{REF:15,REF:17} and spiral search \cite{REF:6}, could be another way to control the intensification step size in the PSS.

Another possible direction includes applying dynamic schemes for the number of population candidates to scale up or down with the global and local search needs. Liang and Juarez (2020) \cite{REF:33} applied a dynamic approach for the population sizing in their algorithm Self-adaptive Virus Optimization Algorithm (SaVOA) and it seems as a promising technique to be applied with the PSS algorithm. Eventually, hybridizing the algorithm with sequential and parallel topologies could be interesting for multimodal optimization problems.

% ********* Reproducibility *********
\section{Reproducibility}
\label{sec:reproducibility}
\par
To promote openness and transparency, the authors have provided the readers with the source code of the algorithm written in C++14, Python 3.7, and Octave (Matlab) programming languages. The source codes can be downloaded from the following online platforms:
\begin{itemize}
    \item Zenodo: \href{https://doi.org/10.5281/zenodo.3630764
}{10.5281/zenodo.3630764}
    \item Github: \href{https://github.com/eesd-epfl/pareto-optimizer}{git@github.com:eesd-epfl/pareto-optimizer .git}
\end{itemize}

% ********* Declarations *********
\section*{Funding}
This work was funded through the Swiss National Science Foundation (SNSF) Project 200021\_175903/1.

\bibliographystyle{spphys}
\bibliography{References}

\end{document}